\documentclass[11pt,feqno]{article}
\usepackage{amssymb,bbm,amsmath,amsthm,amscd}
\usepackage{vmargin}
\usepackage[dvipsnames]{xcolor}
\usepackage{array}
\usepackage{soul}
 
 \usepackage{graphicx}
 \usepackage{subfigure}
\DeclareGraphicsExtensions{.pdf,.eps}

\catcode`\@=11
\newbox\bz@
\newdimen\bdimz@
\def\linethrough#1{\setbox\bz@=\hbox{#1}
\bdimz@=\ht\bz@ \divide\bdimz@ by 5 \advance\bdimz@ by -\dp\bz@ \ht\bz@=\bdimz@
\leavevmode\hbox{$\overline{\overline{\box\bz@}}$\relax}}
\catcode`\@=12

\def\nbZ{{\mathchoice {\hbox{$\sf\textstyle Z\kern-0.4em Z$}}
{\hbox{$\sf\textstyle Z\kern-0.4em Z$}}
{\hbox{$\sf\scriptstyle Z\kern-0.3em Z$}}
{\hbox{$\sf\scriptscriptstyle Z\kern-0.2em Z$}}}}

\newcounter{compteur}[subsection]

\usepackage[colorlinks=true]{hyperref}

\newtheorem{theorem}{Theorem}[section]
\newtheorem{lemma}[theorem]{Lemma}

\newtheorem{corollary}[theorem]{Corollary}

\newtheorem{hypothesis}[theorem]{Hypothesis}
\newtheorem{remark}[theorem]{Remark}

\usepackage{macros}

\makeatletter 

\@addtoreset{equation}{section}

\begin{document}

\author{Mejdi Aza\"iez \thanks{Institut Polytechnique de Bordeaux, Laboratoire I2M CNRS UMR5295, France. {\tt azaiez@enscbp.fr}},
Tom\'as Chac\'on Rebollo \thanks{Departamento EDAN \& IMUS, Universidad de Sevilla, Spain. {\tt chacon@us.es}},
Samuele Rubino \thanks{Departamento EDAN \& IMUS, Universidad de Sevilla, Spain. {\tt samuele@us.es}}}
\title{Streamline derivative projection-based POD-ROM for convection-dominated flows.  Part I : Numerical Analysis}
\maketitle

\begin{abstract}
We introduce improved Reduced Order Models (ROM) for convection-dominated flows. These non-linear closure models are inspired from successful numerical stabilization techniques used in Large Eddy Simulations (LES), such as Local Projection Stabilization (LPS), applied to standard models created by Proper Orthogonal Decomposition (POD) of flows with Galerkin projection. The numerical analysis of the fully Navier--Stokes discretization for the proposed new POD-ROM is presented, by mainly deriving the corresponding error estimates. Also, we suggest an efficient practical implementation of the stabilization term, where the stabilization parameter is approximated by the Discrete Empirical Interpolation Method (DEIM).
\end{abstract}

{\bf{2010 Mathematics Subject Classification:}} Primary 65M12, 65M15, 65M60; \\ Secondary 76D05, 76F20, 76F65.

\medskip

{\bf{Keywords:}} Finite Element Method, Projection Stabilization, Proper Orthogonal Decomposition, Reduced Order Models, Convection-Dominated Flows, Numerical Analysis.

\section{Introduction}\label{sec:Intro}

Reduced-Order Models (ROM) applied to numerical design in modern engineering are a tool that is wide-spreading in the scientific community in the recent years in order to solve complex realistic 
multi-parameters, multi-physics and multi-scale problems, where classical methods such as Finite Difference (FD), Finite Element (FE) or Finite Volume (FV) methods would require up to billions of 
unknowns. On the contrary, ROM are based on a sharp offline/online strategy, and the latter requires a reduced number of unknowns, which allows to face control, optimization, prediction and data analysis problems 
in almost real-time, that is, ultimately, a major goal for industrials. The reduced order modeling offline strategy relies on proper choices for data sampling and
construction of the reduced basis, which will be used then in the online phase, where a proper choice of the reduced model describing the dynamic of the system is
needed. The key feature of ROM is their capability to highly speedup computations, and thus drastically reduce the computational cost of numerical simulations, without compromising too much the
physical accuracy of the solution from the engineering point of view.

\medskip

Among the most popular ROM approaches, Proper Orthogonal Decomposition (POD) strategy provides optimal (from the energetic point of view) basis or modes to represent the dynamics from a given database 
(snapshots) obtained by a full-order system. Onto these reduced basis, a Galerkin projection of the governing equations can be employed to obtain a low-order dynamical system for the basis coefficients. 
The resulting low-order model is named standard POD-ROM, which thus consists in the projection of high-fidelity (full-order) representations of physical problems onto low-dimensional spaces of solutions, 
with a dramatically reduced dimension. These low-dimensional spaces are capable of capturing the dominant characteristics of the solution, their main advantage being that the computations in the 
low-dimensional space can be done at a reduced computational cost. This has led researchers to apply POD-ROM to a variety of physical and engineering problems, including Computational Fluid Dynamics (CFD) 
problems in order to model the Navier--Stokes Equations (NSE), see e.g. \cite{Baiges13, Gunzburger06, Iollo04, Barone12, Beran03, Wang12}. Once applied to the physical problem of interest, POD-ROM can be 
used to solve engineering problems such as shape optimization \cite{Borggaard10a, Borggaard10b} and flow control \cite{Tadmor08, Bergmann08, Peraire99, Tadmor11}.

\medskip

Although POD-ROM can be very computationally efficient and relatively accurate in some flow configurations, they also present several drawbacks. In this report, we address one of them, namely the numerical instability of a straightforward POD-Galerkin procedure applied to convection-dominated flows. The reason of this issue is that, for model reduction purpose, one only keeps few modes that are associated to the large eddies of the flow, which should be sufficient to give a good representation of the kinetic energy of the flow, due to the energetic optimality of the POD basis functions. However, the main amount of viscous dissipation takes place in the small eddies represented by basis functions that are not taken into account, and thus the leading ROM is not able to dissipate enough energy. So, although the disregarded modes do not contain a significant amount of kinetic energy, they have a significant role in the dynamics of the reduced-order system. It is then necessary to close the POD-ROM by modeling the interaction between the computed and the unresolved modes. This problem establishes a parallelism to Large Eddy Simulations (LES) \cite{Sagaut06} of turbulent flows, where the effect of the smallest flow structures on the largest ones is modeled. Since these are also in non-linear interactions, a proper non-linear efficient and accurate closure model should be proposed also in the POD context, considering that in this context the concepts of energy cascade and locality of energy transfer are still valid \cite{Couplet03}.

\medskip

To address this issue, we draw inspiration from the FE context, where stabilized formulations have been developed to deal with the numerical instabilities of the Galerkin method. One of the most popular frameworks for developing stabilized formulations is the Variational Multi-Scale (VMS) method \cite{Hughes98}. In the VMS method, stabilized formulations are obtained by including, in the discretized FE equations, the effect of the part of the solution which cannot be captured by the FE space. This part of the solution
is denoted as the sub-scales. The contribution of the sub-scales turns out not only to be key for the stabilization of the FE problem, but it also allows one to take into account important small-scale effects such as turbulence (cross-stress terms, Reynolds-stress term). Elaborate models for the sub-scales have been developed which allow one to improve the accuracy of VMS stabilized FE methods (cf. \cite{ARCME}). We emphasize that the VMS philosophy is particularly appropriate to the POD setting, in which the hierarchy of small and large structures appears naturally. Indeed, the POD modes are listed in decreasing order of their kinetic energy content.

\medskip

In this work, we propose in particular a POD closure model inspired from successful numerical stabilization techniques used in VMS-LES, such as Local Projection Stabilization (LPS), see \cite{BraackBurman06}. Indeed, the unresolved scales in the proposed stabilized POD-ROM are defined by a projection approach that presents the same structure of the Streamline Derivative-based (SD-based) LPS model (cf. \cite{KnoblochLube09}) in the FE context. This method is an extension to the NSE setting of the one that we proposed, fully analyzed and numerically tested in \cite{ESAIM} for advection-dominated advection-diffusion-reaction equations. Although applications of stabilized methods can already be found in the ROM literature (see \cite{Baiges13, Rozza15, BergmannIollo09, IliescuJohn15, IliescuWang13, IliescuWang14} for the POD context, and also \cite{Maday16, Rozza14} for the Reduced-Basis (RB) context), to the authors' knowledge this is the first time that the SD-based formulation in \cite{KnoblochLube09} has been applied in a POD setting for NSE. The resulting SD-POD-ROM is non-linear, to properly reproduce physical non-linear cross-stress interactions within unresolved and resolved scales, and has a projection-stabilized structure acting only on the high frequencies components of the flow. The structure of the proposed SD-based POD closure model allowed us to perform its numerical analysis for NSE, by mainly deriving error estimates, giving also some hints on how to choose appropriate stabilization para\-meters. In particular, the analy\-sis makes apparent an extra-control on the high frequencies of the convective derivative, aspect of extreme importance, especially when dealing with convection-dominated and turbulent flows. The question of an efficient practical implementation of the strongly non-linear convective stabilization term within the SD-POD-ROM is also addressed, where the non-linear stabilization parameter is approximated using the Discrete Empirical Interpolation Method (DEIM) \cite{Chaturantabut10}. This leads to a computationally efficient and mathematically founded offline/online algorithm (completely separated), implemented over the standard POD-Galerkin ROM. To the best of the authors' knowledge, the use of DEIM for the accurate and efficient computation of the stabilization parameter is new in the literature so far.  

\medskip

The rest of the paper is organized as follows: In Section \ref{sec:VF}, we briefly describe the POD methodology and introduce the new SD-POD-ROM for the NSE. The error analysis for the full discretization (FE in space and backward Euler in time) of the new model is presented in Section \ref{sec:EE}. The practical implementation of the new method is proposed in Section \ref{sec:PI}. Finally, Section \ref{sec:Concl} presents the main conclusions of this work and future research directions.

\section{Streamline derivative projection-based POD-ROM}\label{sec:VF}

We introduce an Initial--Boundary Value Problem (IBVP) for the incompressible evolution Navier--Stokes Equations (NSE). 
For the sake of simplicity, we just impose homogeneous Dirichlet boundary condition on the whole boundary.

\medskip

Let $[0,T]$ be the time interval, and $\Om$ a bounded polyhedral domain in $\mathbb{R}^{d}$, $d=2$ or $3$, with a Lipschitz-continuous boundary $\Ga=\partial\Om$. The transient NSE for an incompressible fluid are given by:

\medskip
\hspace{1cm} {\em Find $\uv:\Om\times (0,T)\longrightarrow\mathbb{R}^{d}$ and $p:\Om\times (0,T)\longrightarrow\mathbb{R}$ such that:}
\BEQ\label{eq:uNS}
\left \{
\begin{array}{rcll}
\partial_{t}\uv + (\uv \cdot\nabla)\uv - \nu \Delta\uv + \nabla p&=&\fv & \qmbx{in} \Om\times (0,T),\\
\div \uv &=&0 & \qmbx{in} \Om\times (0,T),\\
\uv &=& \bf{0} & \qmbx{on} \Ga\times (0,T),\\
\uv(\xv,0)&=&\uv^{0}(\xv) & \qmbx{in} \Om.
\end{array}
\right .
\EEQ
The unknowns are the velocity $\uv(\xv,t)$ and the pressure $p(\xv,t)$ of the incompressible fluid. The data are the source term $\fv(\xv,t)$, which represents a body force per mass unit (typically  the gravity), the kinematic viscosity $\nu$ of the fluid, which is a positive constant, and the initial velocity $\uv^{0}(\xv)$. 

\medskip

To define the weak formulation of problem (\ref{eq:uNS}), we need to introduce some useful notations for functional spaces \cite{Brezis11}. We consider the Sobolev spaces $H^{s}(\Om)$, $s\in \mathbb{R}$, $L^{p}(\Om)$ and $W^{m,p}(\Om)$, $m\in \mathbb{N}$, $1\leq p\leq\infty$. We shall use the following notation for vector-valued Sobolev spaces: ${\bf H}^{s}$, ${\bf L}^{p}$ and ${\bf W}^{m,p}$ respectively shall denote $[H^{s}(\Om)]^{d}$, $[L^{p}(\Om)]^{d}$ and $[W^{m,p}(\Om)]^{d}$ (similarly for tensor spaces of dimension $d\times d$). Also, the parabolic Bochner function spaces $L^{p}(0,T;X)$ and $L^{p}(0,T;{\bf X})$, where $X$ (${\bf X}$) stands for a scalar (vector-valued) Sobolev space shall be denoted by $L^{p}(X)$ and $L^{p}({\bf X})$, respectively. In order to give a variational formulation of problem (\ref{eq:uNS}), let us consider the velocity space:
$$
\Xv={\bf H}_{0}^{1}=[H_{0}^{1}(\Om)]^{d}=\left\{\vv\in [H^{1}(\Om)]^{d} : \vv={\bf 0} \text{ on } \Ga\right\}. 
$$
This is a closed linear subspace of ${\bf H}^{1}$, and thus a Hilbert space endowed with the ${\bf H}^{1}$-norm. 
Thanks to Poincar\'e inequality, the ${\bf H}^{1}$-norm is equivalent on ${\bf H}_{0}^{1}$ to the norm $\nor{\vv}{{\bf H}_{0}^{1}}=\nor{\nabla\vv}{{\bf L}^{2}}$. Also, let us consider the pressure space:
$$
Q=L_{0}^{2}(\Om)=\left\{q\in L^{2}(\Om) : \int_{\Om}q\,d\xv=0\right\}. 
$$
We shall consider the following variational formulation of (\ref{eq:uNS}):

\medskip
\hspace{1cm} {\em Given $\fv\in L^{2}({\bf H}^{-1})$, find $\uv: (0,T)\longrightarrow\Xv$, $p: (0,T)\longrightarrow Q$ such that}
\BEQ\label{eq:fvNS} 
\left \{ 
\begin{array}{rcll}
\disp\frac{d}{dt}(\uv,\vv) + b(\uv,\uv,\vv) + \nu (\nabla\uv,\nabla\vv) - (p,\div\vv) &=&\langle \fv,\vv\rangle & \forall\vv\in\Xv,\quad \text{in } \mathcal{D}'(0,T),\\
(\div \uv,q) &=&0 & \forall q \in Q, \quad \text{a.e. in } (0,T),
\end{array} 
\right .
\EEQ
where $( \cdot,\cdot )$ stands for the $L^2$-inner product in $\Om$, $\langle \cdot,\cdot \rangle$ stands for the duality pairing between $\Xv$ and its dual $\Xv^{'}={\bf H}^{-1}$, and $\mathcal{D}'(0,T)$ is the space of distributions in $(0,T)$. The trilinear form
$b$ is given by: {\it for $\uv,\, \vv,\,\wv \in \Xv$}
\BEQ\label{eq:formb}
b(\uv,\vv,\wv)={1 \over 2} \left[   ( \uv \cdot \nabla\, \vv, \wv) - ( \uv \cdot \nabla\,  \wv, \vv)\right].
\EEQ

\medskip

In order to give a Finite Element (FE) approximation of \eqref{eq:fvNS}, let $\{{\cal T}_{h}\}_{h>0}$ be a family of affine-equivalent, conforming (i.e., without hanging nodes) and regular triangulations of $\overline{\Om}$, 
formed by triangles or quadrilaterals ($d=2$), tetrahedra or hexahedra ($d=3$). For any mesh cell $K \in {\cal T}_{h}$,
its diameter will be denoted by $h_K$ and $h = \max_{K \in {\cal T}_{h}} h_K$. We consider $\Xv^{h}\subset\Xv$, $Q^{h}\subset Q$ being suitable FE spaces for velocity and pressure, respectively. Let us also consider the discrete space of divergence-free functions:
$$
\Vv^{h}=\left\{\vhv\in \Xv^{h} : (\div\vhv,q_{h})=0\quad \forall q_{h}\in Q^{h}\right\}.
$$
The FE approximation of \eqref{eq:fvNS} can be written as follows:
\medskip

\hspace{1cm}{\em Find $\uhv\in\Vv^{h}$ such that}
\BEQ\label{eq:FEapprox}
\disp\frac{d}{dt}(\uhv,\vhv) + b(\uhv,\uhv,\vhv) + \nu (\nabla\uhv,\nabla\vhv) =\langle \fv,\vhv\rangle\quad \forall\vhv\in\Vv^{h}.
\EEQ

\medskip

To ensure error estimates in Theorem \ref{th:PODEE} (main result of the present paper), we have to make the following regularity assumption on the continuous solution $\uv$:
\begin{hypothesis}\label{hp:FERegularity}
In \eqref{eq:fvNS}, assume that $\uv\in L^{\infty}({\bf H}^{2})$, $\partial_{t}\uv\in L^{2}({\bf H}^{1})$ and $\partial_{t}^{2}\uv\in L^{2}({\bf L}^{2})$.
\end{hypothesis}

\subsection{Proper orthogonal decomposition reduced order model}\label{sec:POD-ROM}

We briefly describe the Proper Orthogonal Decomposition (POD) method, following \cite{KunischVolkwein01}. For a detailed presentation, the reader is referred to \cite{Chapelle12, Holmes96, Singler14, Sirovich87, Volkwein11}.

\medskip

Let us consider an ensemble of snapshots $\chi=\text{span}\left\{\uv(\cdot,t_{0}),\ldots,\uv(\cdot,t_{N})\right\}$, which is a collection of velocity data from either numerical simulation results or experimental observations at time $t_{n}=n\Delta t$, $n=0,1,\ldots,N$ and $\Delta t = T/N$. The POD method seeks a low-dimensional basis $\left\{\boldsymbol{\varphi}_{1},\ldots,\boldsymbol{\varphi}_{r}\right\}$ in a real Hilbert space $\mathcal{H}$ that optimally approximates the snapshots in the following sense:
\BEQ\label{eq:PODmeth}
\min\frac{1}{N+1}\sum_{n=0}^{N}\left\| \uv(\cdot,t_{n}) - \sum_{i=1}^{r}\left(\uv(\cdot,t_{n}),\boldsymbol{\varphi}_{i}\right)_{\mathcal{H}}\boldsymbol{\varphi}_{i} \right\|_{\mathcal{H}}^2,
\EEQ
subject to the condition $\left(\boldsymbol{\varphi}_{j},\boldsymbol{\varphi}_{i}\right)_{\mathcal{H}}=\delta_{ij}$, $1\leq i,j \leq r$, where $\delta_{ij}$ is the Kronecker delta. To solve the optimization problem \eqref{eq:PODmeth}, one can consider the eigenvalue problem:
\BEQ\label{eq:eigen}
K\zv_{i}=\lambda_i\zv_{i},\text{ for } 1,\ldots,r,
\EEQ 
where $K\in \mathbb{R}^{(N+1)\times(N+1)}$ is the snapshots correlation matrix with entries: 
$$
K_{mn}=\frac{1}{N+1}\left(\uv(\cdot,t_{n}),\uv(\cdot,t_{m})\right)_{\mathcal{H}},\text{ for } m,n=0,\ldots,N,
$$ 
$\zv_{i}$ is the $i$-th eigenvector, and $\lambda_{i}$ is the associated eigenvalue. The eigenvalues are positive and sorted in descending order $\lambda_{1}\geq\ldots\geq\lambda_{r}>0$. It can be shown that the solution of \eqref{eq:PODmeth}, i.e. the POD basis functions, is given by:
\BEQ\label{eq:PODbasis}
\boldsymbol{\varphi}_{i}(\cdot)=\frac{1}{\sqrt{\lambda_i}}\sum_{n=0}^{N}(\zv_{i})_{n}\uv(\cdot,t_{n}),\quad 1\leq i\leq r,
\EEQ
where $(\zv_{i})_{n}$ is the $n$-th component of the eigenvector $\zv_i$. It can also be shown that the following POD error formula holds \cite{Holmes96, KunischVolkwein01}:
\BEQ\label{eq:PODerr}
\frac{1}{N+1}\sum_{n=0}^{N}\left\| \uv(\cdot,t_{n}) - \sum_{i=1}^{r}\left(\uv(\cdot,t_{n}),\boldsymbol{\varphi}_{i}\right)_{\mathcal{H}}\boldsymbol{\varphi}_{i} \right\|_{\mathcal{H}}^{2} = \sum_{i=r+1}^{M}\lambda_{i},
\EEQ
where $M$ is the rank of $\chi$. Although $\mathcal{H}$ can be any real Hilbert space, in what follows we consider $\mathcal{H}={\bf H}_{0}^{1}$.

\medskip

We consider the following space for the POD setting:
$$
\Xv^{r}=\text{span}\left\{\boldsymbol{\varphi}_{1},\ldots,\boldsymbol{\varphi}_{r}\right\}.
$$

\begin{remark}
Since, as shown in \eqref{eq:PODbasis}, the POD basis functions are linear combinations of the snapshots, the POD basis functions satisfy the boundary conditions in \eqref{eq:uNS} and are solenoidal. If the FE approximations are used as snapshots, the POD basis functions belong to $\Vv^{h}$, which yields $\Xv^{r}\subset\Vv^{h}$.
\end{remark}

The Galerkin projection-based POD-ROM uses both Galerkin truncation and Galerkin projection.
The former yields an approximation of the velocity field by a linear combination of the truncated POD basis:
\BEQ\label{eq:PODsol}
\uv(\xv,t)\approx \uv_{r}(\xv,t)=\sum_{i=1}^{r}a_{i}(t)\boldsymbol{\varphi}_{i}(\xv),
\EEQ
where $\left\{a_{i}(t)\right\}_{i=1}^{r}$
are the sought time-varying coefficients representing the POD-Galerkin trajectories. Note that $r<<\mathcal{N}$, where $\mathcal{N}$ denotes the number of degrees of freedom (d.o.f.) in a full order simulation (e.g., DNS-Direct Numerical Simulation). Replacing the velocity $\uv$ with $\uv_{r}$ in the NSE \eqref{eq:uNS}, using the Galerkin method, and projecting the resulted equations onto the space $\Xv^{r}$, one obtains the standard POD-ROM for the NSE: 
\medskip

\hspace{1cm}{\em Find $\uv_{r}\in\Xv^{r}$ such that} 
\BEQ\label{eq:POD-ROM}
\disp\frac{d}{dt}(\uv_{r},\boldsymbol{\varphi}) + b(\uv_{r},\uv_{r},\boldsymbol{\varphi}) + \nu (\nabla\uv_{r},\nabla\boldsymbol{\varphi}) =\langle \fv,\boldsymbol{\varphi}\rangle\quad \forall\boldsymbol{\varphi}\in\Xv^{r}.
\EEQ
In \eqref{eq:POD-ROM}, the pressure term vanishes due to the fact that all POD modes are solenoidal and satisfy the appropriate boundary conditions. The spatial and temporal discretizations of \eqref{eq:POD-ROM} were considered in \cite{KunischVolkwein02, Luo08}. Despite its appealing computational efficiency, the standard POD-ROM \eqref{eq:POD-ROM} has generally been limited to diffusion-dominated or laminar flows. To overcome this restriction, we develop a non-linear closure model for the standard POD-ROM, which stems from projection-based Variational Multi-Scale (VMS) ideas \cite{IMAJNA, ARCME, CMAME15}.

\subsection{Streamline derivative projection-based method}\label{sec:LPS-POD-ROM}

In projection-based VMS methods, the direct influence of the subgrid-scale model to reproduce the effect of the unresolved scales, usually of (Smagorinsky) eddy viscosity-type in the applications to date \cite{Hughes01a, Hughes01b}, is confined to the small resolved scales. The restriction of the direct influence of the subgrid-scale model to the smaller resolved scales approaches established principles in turbulence theory, namely energy cascade and locality of energy transfer (cf. \cite{Kolmogorov91, Richardson07}). However, for a standard FE discretization, the separation of scales is generally challenging. Indeed, unless special care is taken (e.g., mesh adaptivity is used), the FE basis does not include any a priori information regarding the scales displayed by the underlying problem. On the other hand, note that the hierarchy of basis is implicitly present in a POD setting, since the POD basis functions are already listed in descending order of their kinetic energy content. Thus, the POD represents a perfect setting for the VMS methodology, and VMS closure models for POD-ROM seems to be a natural choice to approximate the effect of the disregarded modes on the retained ones. Indeed, it is well known that a simple Galerkin truncation of POD basis leads to unstable results for convection-dominated and turbulent flows \cite{Aubry88}, and although the disregarded modes do not contain a significant amount of the system's kinetic energy, they have a significant role in the dynamics of the reduced-order system. 

\medskip

To model the effect of the discarded POD modes, various approaches have been proposed, both based on physical insights (cf., e.g., the survey in \cite{Wang12}), or on numerical stabilization techniques for convection-dominated flows (cf. \cite{Baiges13, BergmannIollo09, IliescuJohn15, IliescuWang14}). 
In this paper, we develop an approach that enters in the second group (no ad-hoc eddy viscosity is required, as it is in \cite{Wang12}), and aims to improve the previous works, because on one side a projection-stabilized structure is used (contrary to strategies in \cite{Baiges13, BergmannIollo09, IliescuJohn15}), which allows to act only on the high frequencies components of the flow, and to control them, aspect of extreme importance when dealing especially with convection-dominated and turbulent flows. On the other side, a strongly non-linear closure model is considered here, which is more suitable (with respect to a linear closure model, such as the gradient-based one used in \cite{IliescuWang14}) to reproduce physical non-linear interactions within unresolved and resolved scales. This would allow to improve numerical stability and physical accuracy of the standard Galerkin POD-ROM for convection-dominated flows, with a rather simple driven structure, both for practical implementations such as to perform the numerical analysis. This is not the case, for instance, if we consider a fully residual-based strategy as in \cite{Baiges13, BergmannIollo09}, where the sub-grid terms have a rather complex driven structure, thus increasing computational complexity and setting serious numerical difficulty just to prove stability. The proposed method has been inspired from successful (despite being only weakly consistent) numerical stabilization techniques used in VMS-LES, such as Local Projection Stabilization (LPS), see \cite{BraackBurman06}. Indeed, the unresolved scales in the proposed stabilized POD-ROM are defined by a projection approach that presents the same structure of the Streamline Derivative-based (SD-based) LPS model (cf. \cite{KnoblochLube09}) in the FE context. 

\medskip

To describe our strategy, we define the scalar product:
$$
(\cdot,\cdot)_{\tau}:{\bf L}^{2}\times {\bf L}^{2} \to \mathbb{R},\quad
(\gv,\hv)_{\tau} = \sum_{K\in{\cal T}_{h}}\tau_{K}(\gv,\hv)_{K},
$$
and its associated norm:
$$
\nor{\gv}{\tau}=(\gv,\gv)_{\tau}^{1/2},
$$
where for any $K\in{\cal T}_{h}$, $\tau_{K}$ is a positive local stabilization parameter (to be determined later).
Let us introduce the POD space:
$$
\widehat{\Xv}^{R}=\text{span}\left\{\widehat{\boldsymbol{\varphi}}_{1},\ldots,\widehat{\boldsymbol{\varphi}}_{R}\right\}, \quad R\leq r,
$$
where $\widehat{\boldsymbol{\varphi}}_{i}$, $i=1,\ldots,R$, are the POD basis functions associated to $\widehat{K}$, defined as the snapshots correlation matrix with entries:
\BEQ\label{eq:Khat}
\widehat{K}_{mn}=\frac{1}{N}\left(\uhv^{n+1}\cdot\nabla\uhv^{n+1},\uhv^{m+1}\cdot\nabla\uhv^{m+1}\right),\quad \text{for } m,n=0,\ldots,N-1.
\EEQ
We consider the ${\bf L}^{2}$-orthogonal projection on $\widehat{\Xv}^{R}$, $P_{R}: {\bf L}^{2}\longrightarrow \widehat{\Xv}^{R}$, defined by:
\BEQ\label{eq:PR}
(\uv-P_{R}\uv,\widehat{\vv}_{R})=0,\quad \forall\widehat{\vv}_{R}\in\widehat{\Xv}^{R}.
\EEQ
Let $P_{R}^{\prime}=\mathbb{I}-P_{R}$, where $\mathbb{I}$ is the identity operator. We propose the Streamline Derivative projection-based POD-ROM (SD-POD-ROM) for the NSE: 
\medskip

\hspace{1cm}{\em Find $\uv_{r}\in\Xv^{r}$ such that} 
\BEQ\label{eq:SD-POD-ROM}
\disp\frac{d}{dt}(\uv_{r},\boldsymbol{\varphi}) + b(\uv_{r},\uv_{r},\boldsymbol{\varphi}) + \nu (\nabla\uv_{r},\nabla\boldsymbol{\varphi}) +(P_{R}^{\prime}(\uv_{r}\cdot\nabla\uv_{r}),P_{R}^{\prime}(\uv_{r}\cdot\nabla\boldsymbol{\varphi}))_{\tau}=\langle \fv,\boldsymbol{\varphi}\rangle\quad \forall\boldsymbol{\varphi}\in\Xv^{r}.
\EEQ

\begin{remark}\label{rm:PS-POD-ROM1}
When $\tau_{K}=0$ for any $K\in {\cal T}_{h}$, the SD-POD-ROM \eqref{eq:SD-POD-ROM} coincides with the standard POD-ROM, 
since no numerical dissipation is introduced. When $R = 0$, since numerical diffusion is extended to all the resolved modes $\{\boldsymbol{\varphi}_{1},\ldots,\boldsymbol{\varphi}_{r}\}$,
the SD-POD-ROM \eqref{eq:SD-POD-ROM} becomes a penalty-stabilized method of the form:
\BEQ\label{eq:PS-POD-ROM}
\disp\frac{d}{dt}(\uv_{r},\boldsymbol{\varphi}) + b(\uv_{r},\uv_{r},\boldsymbol{\varphi}) + \nu (\nabla\uv_{r},\nabla\boldsymbol{\varphi}) +(\uv_{r}\cdot\nabla\uv_{r},\uv_{r}\cdot\nabla\boldsymbol{\varphi})_{\tau}=\langle \fv,\boldsymbol{\varphi}\rangle\quad \forall\boldsymbol{\varphi}\in\Xv^{r},
\EEQ
which provides less accuracy with respect to the SD-POD-ROM \eqref{eq:SD-POD-ROM}, see Remark \ref{rm:PS-POD-ROM2} in Section \ref{subsec:MR}. 
\end{remark}

\begin{remark}
Note that the new SD-POD-ROM \eqref{eq:SD-POD-ROM} proposed in the present work is different from the VMS-POD-ROM used in \cite{Wang12}. Indeed, the latter is more specifically based on physical insight, since a sub-grid eddy viscosity of Smagorinsky type is used to model the interactions between the discarded POD modes and those retained in the POD-ROM. Also, it differs from the $P_{R} - VMS-POD-ROM$ introduced in \cite{IliescuWang13}, since there a linear closure model for the standard POD-ROM is considered, which adds artificial viscosity by a term of the form:
$$
\alpha(\overline{P}_{R}^{\prime}(\nabla\uv_{r}),\overline{P}_{R}^{\prime}(\nabla\boldsymbol{\varphi})),
$$  
$\alpha$ being a constant eddy viscosity coefficient, and $\overline{P}_{R}^{\prime}=\mathbb{I}-\overline{P}_{R}$, with $\overline{P}_{R}$ the $L^2$-orthogonal projection on the POD space defined by $\text{span}\{\nabla\boldsymbol{\varphi}_{1},\ldots,\nabla\boldsymbol{\varphi}_{R}\}$. Finally, the proposed SD-POD-ROM \eqref{eq:SD-POD-ROM} is different from the residual-based VMS-POD-ROM introduced in \cite{BergmannIollo09} and the SUPG-POD-ROM introduced in \cite{IliescuJohn15}, since the former uses a projection-stabilized structure, which allows to act only on the high frequencies components of the flow: This guarantees an extra-control on them that prevents high-frequency oscillations without polluting the large scale components of the approximation, see Remark \ref{rm:POD-G-ROM} in Section \ref{subsec:TBpa}. 
\end{remark}

We consider the full discretization of \eqref{eq:SD-POD-ROM}, by using an approximation in time given by the backward Euler method, that is for $n=0,\ldots,N-1$, we compute the approximation $\uv_{r}^{n+1}$ to $\uv^{n+1}=\uv(\cdot,t_{n+1})$ by
\BEQ\label{eq:SD-POD-ROMdisc}
\left \{ 
\begin{array}{lll}
&&\left(\disp\frac{\uv_{r}^{n+1}-\uv_{r}^{n}}{\Delta t},\boldsymbol{\varphi}\right) + b(\uv_{r}^{n+1},\uv_{r}^{n+1},\boldsymbol{\varphi}) + \nu (\nabla\uv_{r}^{n+1},\nabla\boldsymbol{\varphi})\\ \\
&+&(P_{R}^{\prime}(\uv_{r}^{n+1}\cdot\nabla\uv_{r}^{n+1}),P_{R}^{\prime}(\uv_{r}^{n+1}\cdot\nabla\boldsymbol{\varphi}))_{\tau}
=\langle \fv^{n+1},\boldsymbol{\varphi}\rangle\quad \forall\boldsymbol{\varphi}\in\Xv^{r},
\end{array} 
\right .
\EEQ
with $\fv^{n+1}=\fv(\cdot,t_{n+1})$, and the initial condition is given by the elliptic projection of $\uv^{0}$ on $\Xv^{r}$:
\BEQ\label{eq:PODic}
\uv_{r}^{0}=\sum_{i=1}^{r}(\nabla \uv^{0},\nabla\boldsymbol{\varphi}_{i})\boldsymbol{\varphi}_{i}.
\EEQ
In the sequel, we will also denote by $\uhv^{n}$ the FE velocity approximation of \eqref{eq:FEapprox} at $t=t_{n}$.

\medskip

An alternative time discretization could be given by the semi-implicit Euler method:
\BEQ\label{eq:SD-POD-ROMdiscSI}
\left \{ 
\begin{array}{lll}
&&\left(\disp\frac{\uv_{r}^{n+1}-\uv_{r}^{n}}{\Delta t},\boldsymbol{\varphi}\right) + b(\uv_{r}^{n},\uv_{r}^{n+1},\boldsymbol{\varphi}) + \nu (\nabla\uv_{r}^{n+1},\nabla\boldsymbol{\varphi})\\ \\
&+&(P_{R}^{\prime}(\uv_{r}^{n}\cdot\nabla\uv_{r}^{n+1}),P_{R}^{\prime}(\uv_{r}^{n}\cdot\nabla\boldsymbol{\varphi}))_{\tau}
=\langle \fv^{n+1},\boldsymbol{\varphi}\rangle\quad \forall\boldsymbol{\varphi}\in\Xv^{r}.
\end{array} 
\right .
\EEQ
Note that considering a semi-implicit time discretization of the SD-POD-ROM is less costly from the computational point of view with respect to a fully implicit one, which yields a nonlinear algebraic system of equations to be solved. However, the numerical analysis will be performed in detail for the more technical case of the fully implicit time discretization given by \eqref{eq:SD-POD-ROMdisc}.

\section{Error estimates}\label{sec:EE}

In this section, we present the error analysis for the SD-POD-ROM discretization \eqref{eq:SD-POD-ROMdisc}, by mainly focusing on the derivation of error estimates with respect to the continuous solutions $\uv^{n}=\uv(\cdot,t_{n})$, $n=1,\ldots,N$. The error source includes three main components: the spatial FE discretization error, the temporal discretization error, and the POD truncation error. We derive the error estimate in two steps. First, we gather some necessary assumptions and preliminary results in Section \ref{subsec:TBpa}. Then, we present the main result in Section \ref{subsec:MR}.

\subsection{Technical background}\label{subsec:TBpa}

This section provides some technical results that are required for the numerical analysis. Throughout the paper, we shall denote by $C$, $C_1$, $C_2$, $\ldots$ constants that may vary from a line to another, but which are always independent of the FE mesh size $h$, the FE velocity interpolation order $\ell$, the time step $\Delta t$, and the eigenvalues $\lambda_i$.
To prove optimal error estimates in time, we follow \cite{KunischVolkwein01} and include the finite difference quotients $\bar{\partial}\uv^{n}=\disp\frac{\uv^{n}-\uv^{n-1}}{\Delta t}$, for $n=1,\ldots,N$, in the set of snapshots $\chi=\{\uv^{0},\ldots,\uv^{N},\bar{\partial}\uv^{1},\ldots,\bar{\partial}\uv^{N}\}$. As pointed out in \cite{KunischVolkwein01}, the POD error formula \eqref{eq:PODerr} becomes:
\begin{eqnarray}\label{eq:PODerrt}
&&\frac{1}{2N+1}\sum_{n=0}^{N}\left\| \uv^{n} - \sum_{i=1}^{r}\left(\uv^{n},\boldsymbol{\varphi}_{i}\right)_{\mathcal{H}}\boldsymbol{\varphi}_{i} \right\|_{\mathcal{H}}^{2} \nonumber \\
&&+\frac{1}{2N+1}\sum_{n=1}^{N}\left\| \bar{\partial}\uv^{n} - \sum_{i=1}^{r}\left(\bar{\partial}\uv^{n},\boldsymbol{\varphi}_{i}\right)_{\mathcal{H}}\boldsymbol{\varphi}_{i} \right\|_{\mathcal{H}}^{2} 
= \sum_{i=r+1}^{M}\lambda_{i},
\end{eqnarray}
where hereafter $\boldsymbol{\varphi_{i}}$ and $\lambda_{i}$ denote respectively the POD basis functions and the eigenvalues associated to the snapshots correlation matrix with entries: 
$$
K_{m,n}=\frac{1}{2N+1}(\yv^{n},\yv^{m})_{\mathcal{H}},\quad \text{for } m,n=0,\ldots,2N+1,
$$
with $\yv^{i}=\uv^{i}$ for $i=0,\ldots,N+1$, and $\yv^{i}=\bar{\partial}\uv^{i}$ for $i=N+2,\ldots,2N+1$ (but we use the same notation to not overload it).

\medskip

For the subsequent numerical analysis, we need the following technical hypothesis on the stabilization parameters $\tau_{K}$:
\begin{hypothesis}\label{hp:estimCoef}
The stabilization parameters $\tau_{K}$ satisfy the following condition:
\BEQ\label{eq:stabCoef}
\tau_{K}\leq C\,h_{K}^{2},
\EEQ
for all $K\in{\cal T}_{h}$, and a positive constant $C$ independent of $h$.
\end{hypothesis} 

\begin{remark}
The question whether the stabilization parameters should depend on the spatial resolution of the underlying FE space, or on the number of POD basis functions used has been addressed in \cite{IliescuJohn15}, by means of numerical analysis arguments. In that work, numerical investigations using both definitions suggested that the one based on estimates from the underlying FE discretization provides a better suppression of numerical oscillations, and thus guarantees a more effective numerical stabilization. For this reason, we make here assumption \ref{hp:estimCoef} on the stabilization parameters, which is also essential for the subsequent numerical analysis.
\end{remark}

\begin{hypothesis}\label{hp:FEerrEst}
Assume that the FE approximation $\uhv^{n}$ of \eqref{eq:FEapprox} satisfies the following error estimate:
\BEQ\label{eq:FEvelEE}
\nor{\uv-\uhv}{l^{\infty}({\bf L}^{2})}+\nor{\nabla(\uv-\uhv)}{l^{2}({\bf L}^{2})}+\nor{\uv\cdot\nabla\uv-\uhv\cdot\nabla\uhv}{l^{2}(\tau)}\leq C(h^{\ell}+\Delta t).
\EEQ
Assume also the following standard approximation property (see, e.g., page 166 in \cite{Layton08}):
\BEQ\label{eq:FEpresEE}
\inf_{q_{h}\in Q^{h}}\nor{p-q_{h}}{{\bf L}^{2}}\leq C\,h^{\ell}.
\EEQ
\end{hypothesis}
\begin{lemma}\label{lm:stabPR}
Assume that Hypothesis \ref{hp:estimCoef} holds. Then, for all $\gv\in {\bf L}^{2}$, the following estimate is satisfied:
\BEQ\label{eq:stabBound}
\nor{P_{R}^{\prime}(\gv)}{\tau}\leq C\,h\nor{\gv}{{\bf L}^{2}}.
\EEQ
\end{lemma}
{\bf Proof.}
By using \eqref{eq:stabCoef} and the stability of $P_{R}$ in the ${\bf L}^{2}$-norm, it follows:
$$
\nor{P_{R}^{\prime}(\gv)}{\tau}^{2}\leq C\,h^{2}\nor{P_{R}^{\prime}(\gv)}{{\bf L}^{2}}^{2}\leq C\,h^{2}\nor{\gv}{{\bf L}^{2}}^{2}.
$$
Thus, the estimate (\ref{eq:stabBound}) can be deduced.
\qed

\medskip

We have the following error estimate for $\vv_{r}^{n}=\sum_{i=1}^{r}(\nabla\uv^{n},\nabla\boldsymbol{\varphi}_{i})\boldsymbol{\varphi}_{i}$, i.e. the elliptic projection of $\uv^{n}$ on $\Xv^{r}$ (see \cite{IliescuWang14}, Lemma 3.3):
\begin{lemma}\label{lm:PODlm}
\BEQ\label{eq:PODlmL2}
\frac{1}{N+1}\sum_{n=0}^{N}\nor{\uv^{n}-\vv_{r}^{n}}{{\bf L}^{2}}^{2}\leq C\left(h^{2\ell} + \Delta t ^{2} + \sum_{i=r+1}^{M}\lambda_{i}\right),
\EEQ
\BEQ\label{eq:PODlmH01}
\frac{1}{N+1}\sum_{n=0}^{N}\nor{\nabla(\uv^{n}-\vv_{r}^{n})}{{\bf L}^{2}}^{2}\leq C\left(h^{2\ell} + \Delta t ^{2} + \sum_{i=r+1}^{M}\lambda_{i}\right).
\EEQ
\end{lemma}

\begin{corollary}\label{co:PODco}
\BEQ\label{eq:PODlmL2t}
\frac{1}{N}\sum_{n=1}^{N}\nor{\partial_{t}(\uv^{n}-\vv_{r}^{n})}{{\bf L}^{2}}^{2}\leq C\left(h^{2\ell} + \Delta t ^{2} + \sum_{i=r+1}^{M}\lambda_{i}\right).
\EEQ
\end{corollary}
The proof of this corollary follows along the same lines as the proof of Lemma \ref{lm:PODlm}. Note that
it is exactly at this point that we use the fact that the finite difference quotients $\bar{\partial}\uv^{n}$ are included in the set of snapshots (see Remark 1 in \cite{KunischVolkwein01}).
\begin{lemma}[See Lemma 13 in \cite{Layton08}]\label{lm:formb}
For any function $\uv,\vv,\wv\in\Xv$, the skew-symmetric trilinear form $b(\cdot,\cdot,\cdot)$ satisfies:
\BEQ\label{eq:formb}
b(\uv,\vv,\vv)=0,
\EEQ 
\BEQ\label{eq:formb1}
b(\uv,\vv,\wv)\leq C\nor{\nabla\uv}{{\bf L}^{2}}\nor{\nabla\vv}{{\bf L}^{2}}\nor{\nabla\wv}{{\bf L}^{2}}.
\EEQ
\end{lemma}

We have the following existence and stability result for the SD-POD-ROM \eqref{eq:SD-POD-ROMdisc}:
\begin{lemma}\label{lm:stabEst}
Problem \eqref{eq:SD-POD-ROMdisc} admits a solution that satisfies the following bound:
\BEQ\label{eq:stabEst}
\nor{\uv_{r}^{k}}{{\bf L}^{2}}^{2}+\Delta t\sum_{n=0}^{N-1}\left(\nu\nor{\nabla\uv_{r}^{n+1}}{{\bf L}^{2}}^{2} 
+ \nor{P_{R}^{\prime}(\uv_{r}^{n+1}\cdot\nabla\uv_{r}^{n+1})}{\tau}^{2}\right)\leq
\nor{\uv_{r}^{0}}{{\bf L}^{2}}^{2}+\frac{\Delta t}{\nu}\sum_{n=0}^{N-1}\nor{\fv^{n+1}}{{\bf H}^{-1}}^{2},
\EEQ
for $k=0,\ldots,N$.
\end{lemma}
{\bf Proof.}
Problem \eqref{eq:SD-POD-ROMdisc} can be written as:
\BEQ\label{eq:SD-POD-ROMdiscR}
b(\uv_{r}^{n+1},\uv_{r}^{n+1},\boldsymbol{\varphi})+\widetilde{a}(\uv_{r}^{n+1},\boldsymbol{\varphi})+
(P_{R}^{\prime}(\uv_{r}^{n+1}\cdot\nabla\uv_{r}^{n+1}),P_{R}^{\prime}(\uv_{r}^{n+1}\cdot\nabla\boldsymbol{\varphi}))_{\tau}=
\langle \widetilde{\fv}^{n+1},\boldsymbol{\varphi}\rangle\quad \forall\boldsymbol{\varphi}\in\Xv^{r},
\EEQ
where $\widetilde{a}(\uv_{r}^{n+1},\boldsymbol{\varphi})=\Delta t^{-1}(\uv_{r}^{n+1},\boldsymbol{\varphi})+\nu (\nabla\uv_{r}^{n+1},\nabla\boldsymbol{\varphi})$, and $ \widetilde{\fv}^{n+1}=\fv^{n+1}+\Delta t^{-1}(\uv_{r}^{n},\boldsymbol{\varphi})$. 
This problem fits into the same functional framework as for implicit discretizations of the steady NSE (cf. \cite{CMAME15} for instance), since $\widetilde{a}$ is an inner product on space $\Xv$ that generates a norm equivalent to the ${\bf H}^{1}$-norm.
Then, the existence of a solution follows from Brouwer's fixed point theorem \cite{Brezis11} (see Steps 1 and 2 of Theorem 3.6 in \cite{CMAME15} for instance).

\medskip

To prove estimate \eqref{eq:stabEst}, we choose $\boldsymbol{\varphi}=\uv_{r}^{n+1}$ in \eqref{eq:SD-POD-ROMdisc}, and note $b(\uv_{r}^{n+1},\uv_{r}^{n+1},\uv_{r}^{n+1})=0$ by \eqref{eq:formb}, so that we obtain:
\BEQ\label{eq:stabEst1}
\left \{ 
\begin{array}{lll}
&&(\uv_{r}^{n+1}-\uv_{r}^{n},\uv_{r}^{n+1}) + \nu\Delta t (\nabla\uv_{r}^{n+1},\nabla\uv_{r}^{n+1}) \\ \\
&+&\Delta t(P_{R}^{\prime}(\uv_{r}^{n+1}\cdot\nabla\uv_{r}^{n+1}),P_{R}^{\prime}(\uv_{r}^{n+1}\cdot\nabla\uv_{r}^{n+1}))_{\tau}
=\Delta t\langle \fv^{n+1},\uv_{r}^{n+1}\rangle.
\end{array} 
\right .
\EEQ
Using the identity:
$$
(a-b)a=\frac{1}{2}(|a|^{2}-|b|^{2}+|a-b|^{2}),\quad\forall a,b\in \mathbb{R},
$$
and Young's inequality, from \eqref{eq:stabEst1} we get:
\BEQ\label{eq:stabEst2}
\nor{\uv_{r}^{n+1}}{{\bf L}^{2}}^{2}-\nor{\uv_{r}^{n}}{{\bf L}^{2}}^{2}+\nu\Delta t\nor{\nabla\uv_{r}^{n+1}}{{\bf L}^{2}}^{2}+\Delta t\nor{P_{R}^{\prime}(\uv_{r}^{n+1}\cdot\nabla\uv_{r}^{n+1})}{\tau}^{2}
\leq \disp\frac{\Delta t}{\nu}\nor{\fv^{n+1}}{{\bf H}^{-1}}^{2}.
\EEQ
Then, the stability estimate \eqref{eq:stabEst} follows by summing \eqref{eq:stabEst2} from $n=0$ to $k\leq N-1$.
\qed

\begin{remark}\label{rm:POD-G-ROM}
The stability estimate \eqref{eq:stabEst}, which makes apparent the estimate of the
convective stabilization term, guarantees an extra-control on the high frequencies of the convective
derivative, which is not obtained by the standard Galerkin POD-ROM. This is an aspect of extreme importance, especially when dealing with convection-dominated flows.
\end{remark}

\subsection{Error estimate  for the SD-POD-ROM}\label{subsec:MR}

We are now in position to prove the following error estimate result for the SD-POD-ROM defined by  \eqref{eq:SD-POD-ROMdisc}:
\begin{theorem}\label{th:PODEE}
Under the regularity assumption on the continuous solution (Hypothesis \ref{hp:FERegularity}), the assumption on the FE approximation (Hypothesis \ref{hp:FEerrEst}), the assumption on 
the stabilization parameters (Hypothesis \ref{hp:estimCoef}), and supposing that $\nor{\uv^{0}-\uv_{r}^{0}}{{\bf L}^{2}}=\mathcal{O}(h^{\ell})$, the solution of the SD-POD-ROM \eqref{eq:SD-POD-ROMdisc}
satisfies the following error estimate: For a sufficiently large $r:\lambda_{r+1}\sim\mathcal{O}(\nu)$, there exists $\Delta t^{*}>0$ such that the inequality
\begin{eqnarray}\label{eq:PODEE}
&&\frac{1}{N+1}\sum_{n=0}^{N}\nor{\uv^{n}-\uv_{r}^{n}}{{\bf L}^{2}}^{2}+\nu\Delta t\sum_{n=0}^{N-1}\nor{\nabla(\uv^{n+1}-\uv_{r}^{n+1})}{{\bf L}^{2}}^{2}\nonumber \\
&\leq& C\left(h^{2\ell}+\Delta t^{2} + \sum_{i=r+1}^{M}\lambda_{i}+h^{2}\sum_{i=R+1}^{M}\widehat{\lambda}_{i}\right),
\end{eqnarray}
holds for $h\sim\Delta t \left(\sim\mathcal{O}(\sqrt{\nu})\right) \leq \Delta t^{*}$, where $\Delta t^{*}$ will be determined throughout the proof, and $\widehat{\lambda}_{i}$, $i=R+1,\ldots,M$, in the right-hand side of \eqref{eq:PODEE} are the eigenvalues associated to the snapshots correlation matrix $\widehat{K}$ previously defined in \eqref{eq:Khat}.
\end{theorem}

{\bf Proof.} 
We start deriving the error bound by splitting the error into two terms:
\BEQ\label{eq:splitErr}
\uv^{n+1}-\uv_{r}^{n+1}=(\uv^{n+1}-\vv_{r}^{n+1})-(\uv_{r}^{n+1}-\vv_{r}^{n+1})=\boldsymbol{\eta}^{n+1}-\boldsymbol{\phi}_{r}^{n+1}.
\EEQ
The first term, $\boldsymbol{\eta}^{n+1}=\uv^{n+1}-\vv_{r}^{n+1}$, represents the difference between $\uv^{n+1}$ and its elliptic projection on $\Xv^r$. The second term, $\boldsymbol{\phi}_{r}^{n+1}=\uv_{r}^{n+1}-\vv_{r}^{n+1}$, is the remainder.

\medskip

Next, we construct the error equation. We first evaluate the weak formulation of the NSE \eqref{eq:fvNS} at $t=t_{n+1}$, and let $\vv=\boldsymbol{\varphi}_{r}$, then subtract the 
SD-POD-ROM \eqref{eq:SD-POD-ROMdisc} from it. For any $\boldsymbol{\varphi}_{r}\in\Xv^r$, we obtain:
\begin{eqnarray}\label{eq:ErrEq}
&&(\partial_{t}\uv^{n+1},\boldsymbol{\varphi}_{r})-\frac{1}{\Delta t}(\uv_{r}^{n+1}-\uv_{r}^{n},\boldsymbol{\varphi}_{r})+
\nu\left(\nabla(\uv^{n+1}-\uv_{r}^{n+1}),\nabla\boldsymbol{\varphi}_{r}\right)\nonumber \\ \nonumber \\
&+&b(\uv^{n+1},\uv^{n+1},\boldsymbol{\varphi}_{r})-b(\uv_{r}^{n+1},\uv_{r}^{n+1},\boldsymbol{\varphi}_{r})-(p^{n+1},\div\boldsymbol{\varphi}_{r})\nonumber \\ \nonumber \\
&-&\left(P_{R}^{\prime}(\uv_{r}^{n+1}\cdot\nabla\uv_{r}^{n+1}),P_{R}^{\prime}(\uv_{r}^{n+1}\cdot\nabla\boldsymbol{\varphi}_{r})\right)_{\tau}=0.
\end{eqnarray}
By adding and subtracting the different quotient term $\frac{1}{\Delta t}(\uv^{n+1}-\uv^{n},\boldsymbol{\varphi}_{r})$ in \eqref{eq:ErrEq}, and applying the decomposition \eqref{eq:splitErr}, 
we get, for any $\boldsymbol{\varphi}_{r}\in\Xv^r$:
\begin{eqnarray}\label{eq:ErrEq0}
&&(\partial_{t}\uv^{n+1}-\frac{\uv^{n+1}-\uv^{n}}{\Delta t},\boldsymbol{\varphi}_{r})+\frac{1}{\Delta t}(\boldsymbol{\eta}^{n+1}-\boldsymbol{\phi}_{r}^{n+1},\boldsymbol{\varphi}_{r})
-\frac{1}{\Delta t}(\boldsymbol{\eta}^{n}-\boldsymbol{\phi}_{r}^{n},\boldsymbol{\varphi}_{r})\nonumber \\ \nonumber \\
&+&\nu\left(\nabla(\boldsymbol{\eta}^{n+1}-\boldsymbol{\phi}_{r}^{n+1}),\nabla\boldsymbol{\varphi}_{r}\right)+b(\uv^{n+1},\uv^{n+1},\boldsymbol{\varphi}_{r})-b(\uv_{r}^{n+1},\uv_{r}^{n+1},\boldsymbol{\varphi}_{r})-(p^{n+1},\div\boldsymbol{\varphi}_{r})\nonumber \\ \nonumber \\
& &-\left(P_{R}^{\prime}(\uv_{r}^{n+1}\cdot\nabla\uv_{r}^{n+1}),P_{R}^{\prime}(\uv_{r}^{n+1}\cdot\nabla\boldsymbol{\varphi}_{r})\right)_{\tau}=0.
\end{eqnarray}
Note that $(\nabla\boldsymbol{\eta}^{n+1},\nabla\boldsymbol{\varphi}_{r})=0$, since $\vv_{r}^{n+1}$ is the elliptic projection of $\uv^{n+1}$ on $\Xv^r$. 
Choosing $\boldsymbol{\varphi}_{r}=\boldsymbol{\phi}_{r}^{n+1}$ in \eqref{eq:ErrEq0} and letting $\boldsymbol{r}^{n}=\partial_{t}\uv^{n+1}-\frac{\uv^{n+1}-\uv^{n}}{\Delta t}$, we obtain:
\begin{eqnarray}\label{eq:ErrEq1}
&&\frac{1}{\Delta t}(\boldsymbol{\phi}_{r}^{n+1}-\boldsymbol{\phi}_{r}^{n},\boldsymbol{\phi}_{r}^{n+1})
+\nu(\nabla\boldsymbol{\phi}_{r}^{n+1},\nabla\boldsymbol{\phi}_{r}^{n+1})\nonumber \\ \nonumber \\
&=&\langle\boldsymbol{r}^{n},\boldsymbol{\phi}_{r}^{n+1}\rangle
+\frac{1}{\Delta t}(\boldsymbol{\eta}^{n+1}-\boldsymbol{\eta}^{n},\boldsymbol{\phi}_{r}^{n+1})
+b(\uv^{n+1},\uv^{n+1},\boldsymbol{\phi}_{r}^{n+1})-b(\uv_{r}^{n+1},\uv_{r}^{n+1},\boldsymbol{\phi}_{r}^{n+1})
\nonumber \\ \nonumber \\
&&-(p^{n+1},\div\boldsymbol{\phi}_{r}^{n+1})-\left(P_{R}^{\prime}(\uv_{r}^{n+1}\cdot\nabla\uv_{r}^{n+1}),P_{R}^{\prime}(\uv_{r}^{n+1}\cdot\nabla\boldsymbol{\phi}_{r}^{n+1})\right)_{\tau}.
\end{eqnarray}
First, we estimate the left-hand side of \eqref{eq:ErrEq1}, by applying Cauchy--Schwarz and Young's inequalities:
\BEQ\label{eq:ErrEqLHS}
\frac{1}{\Delta t}\nor{\boldsymbol{\phi}_{r}^{n+1}}{{\bf L}^{2}}^{2}-\frac{1}{\Delta t}(\boldsymbol{\phi}_{r}^{n},\boldsymbol{\phi}_{r}^{n+1})
+\nu\nor{\nabla\boldsymbol{\phi}_{r}^{n+1}}{{\bf L}^{2}}^{2}\geq
\frac{1}{2\Delta t}(\nor{\boldsymbol{\phi}_{r}^{n+1}}{{\bf L}^{2}}^{2}-\nor{\boldsymbol{\phi}_{r}^{n}}{{\bf L}^{2}}^{2})
+\nu\nor{\nabla\boldsymbol{\phi}_{r}^{n+1}}{{\bf L}^{2}}^{2}.
\EEQ
Next, we estimate the terms on the right-hand side of \eqref{eq:ErrEq1} one by one. Using Cauchy--Schwarz and Young's inequalities, we get for the first two terms on the right-hand side of \eqref{eq:ErrEq1}:
\BEQ\label{eq:timederErrEq}
\langle\boldsymbol{r}^{n},\boldsymbol{\phi}_{r}^{n+1}\rangle\leq \nor{\boldsymbol{r}^{n}}{{\bf H}^{-1}}\nor{\nabla\boldsymbol{\phi}_{r}^{n+1}}{{\bf L}^2}\leq \frac{\varepsilon^{-1}}{4}\nor{\boldsymbol{r}^{n}}{{\bf H}^{-1}}^{2}+\varepsilon\nor{\nabla\boldsymbol{\phi}_{r}^{n+1}}{{\bf L}^{2}}^{2},
\EEQ
\begin{eqnarray}\label{eq:diffErrEq}
\frac{1}{\Delta t}(\boldsymbol{\eta}^{n+1}-\boldsymbol{\eta}^{n},\boldsymbol{\phi}_{r}^{n+1})&\leq& C_{P}\left\| \frac{1}{\Delta t}(\boldsymbol{\eta}^{n+1}-\boldsymbol{\eta}^{n})\right\|_{{\bf L}^2}\nor{\nabla\boldsymbol{\phi}_{r}^{n+1}}{{\bf L}^2}\nonumber \\
&\leq& \frac{\varepsilon^{-1}C_{P}^{2}}{4}\left\| \frac{1}{\Delta t}(\boldsymbol{\eta}^{n+1}-\boldsymbol{\eta}^{n})\right\|_{{\bf L}^2}^{2}+\varepsilon\nor{\nabla\boldsymbol{\phi}_{r}^{n+1}}{{\bf L}^{2}}^{2}, 
\end{eqnarray}
for some small positive constant $\varepsilon$, and $C_{P}$ denoting the Poincar\'e constant. 

\medskip

The nonlinear convective terms in \eqref{eq:ErrEq1} can be written as follows:
\begin{eqnarray}\label{eq:convErrEq}
&&b(\uv^{n+1},\uv^{n+1},\boldsymbol{\phi}_{r}^{n+1})-b(\uv_{r}^{n+1},\uv_{r}^{n+1},\boldsymbol{\phi}_{r}^{n+1})\nonumber \\
&=&b(\uv_{r}^{n+1},\boldsymbol{\eta}^{n+1}-\boldsymbol{\phi}_{r}^{n+1},\boldsymbol{\phi}_{r}^{n+1})
+b(\boldsymbol{\eta}^{n+1}-\boldsymbol{\phi}_{r}^{n+1},\uv^{n+1},\boldsymbol{\phi}_{r}^{n+1})\nonumber \\
&=&b(\uv_{r}^{n+1},\boldsymbol{\eta}^{n+1},\boldsymbol{\phi}_{r}^{n+1})
+b(\boldsymbol{\eta}^{n+1},\uv^{n+1},\boldsymbol{\phi}_{r}^{n+1})
-b(\boldsymbol{\phi}_{r}^{n+1},\uv^{n+1},\boldsymbol{\phi}_{r}^{n+1}),
\end{eqnarray}
where we have used $b(\uv_{r}^{n+1},\boldsymbol{\phi}_{r}^{n+1},\boldsymbol{\phi}_{r}^{n+1})=0$, which follows from \eqref{eq:formb}. Next, we estimate each term on the right-hand side 
of \eqref{eq:convErrEq}. Since $\uv_{r}^{n+1},\boldsymbol{\eta}^{n+1},\boldsymbol{\phi}_{r}^{n+1},\uv^{n+1}\in\Xv$, 
we can apply the standard bound \eqref{eq:formb1} for the trilinear form $b(\cdot,\cdot,\cdot)$, and use Young's inequality to get:
\begin{eqnarray}\label{eq:conv1}
b(\uv_{r}^{n+1},\boldsymbol{\eta}^{n+1},\boldsymbol{\phi}_{r}^{n+1})&\leq&
C\nor{\nabla\uv_{r}^{n+1}}{{\bf L}^{2}}\nor{\nabla\boldsymbol{\eta}^{n+1}}{{\bf L}^{2}}\nor{\nabla\boldsymbol{\phi}_{r}^{n+1}}{{\bf L}^{2}}\nonumber \\ 
&\leq& \frac{\varepsilon^{-1}C^{2}}{4}\nor{\nabla\uv_{r}^{n+1}}{{\bf L}^{2}}^{2}\nor{\nabla\boldsymbol{\eta}^{n+1}}{{\bf L}^{2}}^{2} +\varepsilon\nor{\nabla\boldsymbol{\phi}_{r}^{n+1}}{{\bf L}^{2}}^{2};
\end{eqnarray}
\begin{eqnarray}\label{eq:conv2}
b(\boldsymbol{\eta}^{n+1},\uv^{n+1},\boldsymbol{\phi}_{r}^{n+1})&\leq&
C\nor{\nabla\boldsymbol{\eta}^{n+1}}{{\bf L}^{2}}\nor{\nabla\uv^{n+1}}{{\bf L}^{2}}\nor{\nabla\boldsymbol{\phi}_{r}^{n+1}}{{\bf L}^{2}}\nonumber \\ 
&\leq& \frac{\varepsilon^{-1}C^{2}}{4}\nor{\nabla\uv^{n+1}}{{\bf L}^{2}}^{2}\nor{\nabla\boldsymbol{\eta}^{n+1}}{{\bf L}^{2}}^{2} +\varepsilon\nor{\nabla\boldsymbol{\phi}_{r}^{n+1}}{{\bf L}^{2}}^{2}.
\end{eqnarray}
For the last nonlinear convective term, applying H\"older's inequality, Sobolev embedding theorem, and Young's inequality yields:
\begin{eqnarray}\label{eq:conv3}
b(\boldsymbol{\phi}_{r}^{n+1},\uv^{n+1},\boldsymbol{\phi}_{r}^{n+1})&\leq&
C\nor{\boldsymbol{\phi}_{r}^{n+1}}{{\bf L}^{2}}(\nor{\nabla\uv^{n+1}}{{\bf L}^{3}}+\nor{\uv^{n+1}}{{\bf L}^{\infty}})\nor{\nabla\boldsymbol{\phi}_{r}^{n+1}}{{\bf L}^{2}}\nonumber \\  
&\leq& \frac{\varepsilon^{-1}C^{2}}{4}\nor{\nabla\uv^{n+1}}{{\bf H}^{1}}^{2}\nor{\boldsymbol{\phi}_{r}^{n+1}}{{\bf L}^{2}}^{2} 
+\varepsilon\nor{\nabla\boldsymbol{\phi}_{r}^{n+1}}{{\bf L}^{2}}^{2}.
\end{eqnarray}

Since $\boldsymbol{\phi}_{r}^{n+1}\in \Xv^{r}\subset \Vv^{h}$, the pressure term on the right-hand side of \eqref{eq:ErrEq1} can be written as:
$$
-(p^{n+1},\div\boldsymbol{\phi}_{r}^{n+1})=-(p^{n+1}-q_{h},\div\boldsymbol{\phi}_{r}^{n+1}),
$$
for any $q_{h}\in Q^{h}$. Thus, the pressure term can be estimated as follows, by using Cauchy--Schwarz and Young's inequalities:
\BEQ\label{eq:presErrEq}
-(p^{n+1},\div\boldsymbol{\phi}_{r}^{n+1})\leq \frac{\varepsilon^{-1}}{4}\nor{p^{n+1}-q_{h}}{{\bf L}^{2}}^{2}
+\varepsilon\nor{\nabla\boldsymbol{\phi}_{r}^{n+1}}{{\bf L}^{2}}^{2}.
\EEQ

The last term on the right-hand side of \eqref{eq:ErrEq1} can be estimated   using Cauchy--Schwarz and Young's inequalities:
\begin{eqnarray}\label{eq:stabErrEq}
&-&\left(P_{R}^{\prime}(\uv_{r}^{n+1}\cdot\nabla\uv_{r}^{n+1}),P_{R}^{\prime}(\uv_{r}^{n+1}\cdot\nabla\boldsymbol{\phi}_{r}^{n+1})\right)_{\tau} \nonumber \\  
&=&\left(P_{R}^{\prime}(\uv_{r}^{n+1}\cdot\nabla\boldsymbol{\eta}^{n+1}),P_{R}^{\prime}(\uv_{r}^{n+1}\cdot\nabla\boldsymbol{\phi}_{r}^{n+1})\right)_{\tau}
-\left(P_{R}^{\prime}(\uv_{r}^{n+1}\cdot\nabla\boldsymbol{\phi}_{r}^{n+1}),P_{R}^{\prime}(\uv_{r}^{n+1}\cdot\nabla\boldsymbol{\phi}_{r}^{n+1})\right)_{\tau}\nonumber \\
& &-\left(P_{R}^{\prime}(\uv_{r}^{n+1}\cdot\nabla\uv^{n+1}),P_{R}^{\prime}(\uv_{r}^{n+1}\cdot\nabla\boldsymbol{\phi}_{r}^{n+1})\right)_{\tau}\nonumber \\ 
&\leq&\nor{P_{R}^{\prime}(\uv_{r}^{n+1}\cdot\nabla\boldsymbol{\eta}^{n+1})}{\tau}\nor{P_{R}^{\prime}(\uv_{r}^{n+1}\cdot\nabla\boldsymbol{\phi}_{r}^{n+1})}{\tau}
-\nor{P_{R}^{\prime}(\uv_{r}^{n+1}\cdot\nabla\boldsymbol{\phi}_{r}^{n+1})}{\tau}^{2}\nonumber \\
&&+\nor{P_{R}^{\prime}(\uv_{r}^{n+1}\cdot\nabla\uv^{n+1})}{\tau}\nor{P_{R}^{\prime}(\uv_{r}^{n+1}\cdot\nabla\boldsymbol{\phi}_{r}^{n+1})}{\tau}\nonumber \\ 
&\leq& \nor{P_{R}^{\prime}(\uv_{r}^{n+1}\cdot\nabla\boldsymbol{\eta}^{n+1})}{\tau}^{2}-\frac{1}{2}\nor{P_{R}^{\prime}(\uv_{r}^{n+1}\cdot\nabla\boldsymbol{\phi}_{r}^{n+1})}{\tau}^{2}
+\nor{P_{R}^{\prime}(\uv_{r}^{n+1}\cdot\nabla\uv^{n+1})}{\tau}^{2}.
\end{eqnarray}
Substituting inequalities \eqref{eq:timederErrEq}-\eqref{eq:diffErrEq} and \eqref{eq:conv1}-\eqref{eq:stabErrEq} in \eqref{eq:ErrEq1}, multiplying by $2\Delta t$ both sides and taking $\varepsilon=\nu/6$, we obtain:
\begin{eqnarray}\label{eq:ErrEq3}
&&\nor{\boldsymbol{\phi}_{r}^{n+1}}{{\bf L}^{2}}^{2}-\nor{\boldsymbol{\phi}_{r}^{n}}{{\bf L}^{2}}^{2}
+\nu\Delta t\nor{\nabla\boldsymbol{\phi}_{r}^{n+1}}{{\bf L}^{2}}^{2}
+\Delta t\nor{P_{R}^{\prime}(\uv_{r}^{n+1}\cdot\nabla\boldsymbol{\phi}_{r}^{n+1})}{\tau}^{2}\nonumber \\
&\leq&\frac{3}{\nu}\,\Delta t\nor{\boldsymbol{r}^{n}}{{\bf H}^{-1}}^{2}
+\frac{3}{\nu}C_{P}^{2}\,\Delta t\left\| \frac{1}{\Delta t}(\boldsymbol{\eta}^{n+1}-\boldsymbol{\eta}^{n})\right\|_{{\bf L}^2}^{2}
+\frac{3\,C^{2}}{\nu}\Delta t\nor{\nabla\uv_{r}^{n+1}}{{\bf L}^{2}}^{2}\nor{\nabla\boldsymbol{\eta}^{n+1}}{{\bf L}^{2}}^{2}
\nonumber \\ 
&&+\frac{3\,C^{2}}{\nu}\Delta t\nor{\nabla\uv^{n+1}}{{\bf L}^{2}}^{2}\nor{\nabla\boldsymbol{\eta}^{n+1}}{{\bf L}^{2}}^{2}
+\frac{3\,C^{2}}{\nu}\Delta t\nor{\nabla\uv^{n+1}}{{\bf H}^{1}}^{2}\nor{\boldsymbol{\phi}_{r}^{n+1}}{{\bf L}^{2}}^{2}
\nonumber \\
&&+\frac{3}{\nu}\Delta t\nor{p^{n+1}-q_{h}}{{\bf L}^{2}}^{2}
+2\,\Delta t\left(\nor{P_{R}^{\prime}(\uv_{r}^{n+1}\cdot\nabla\boldsymbol{\eta}^{n+1})}{\tau}^{2}+\nor{P_{R}^{\prime}(\uv_{r}^{n+1}\cdot\nabla\uv^{n+1})}{\tau}^{2}\right)\nonumber \\
&\leq& C_{1}\,\Delta t\nor{\boldsymbol{r}^{n}}{{\bf H}^{-1}}^{2}
+C_{1}\,\Delta t\left\| \frac{1}{\Delta t}(\boldsymbol{\eta}^{n+1}-\boldsymbol{\eta}^{n})\right\|_{{\bf L}^2}^{2}\nonumber \\
&&+C_{1}\,\Delta t\left(\nor{\nabla\uv_{r}^{n+1}}{{\bf L}^{2}}^{2}
+\nor{\nabla\uv^{n+1}}{{\bf L}^{2}}^{2}\right)\nor{\nabla\boldsymbol{\eta}^{n+1}}{{\bf L}^{2}}^{2}\nonumber \\
&&+C_{1}\,\Delta t\nor{\nabla\uv^{n+1}}{{\bf H}^{1}}^{2}\nor{\boldsymbol{\phi}_{r}^{n+1}}{{\bf L}^{2}}^{2}
+C_{1}\,\Delta t \nor{p^{n+1}-q_{h}}{{\bf L}^{2}}^{2} \nonumber \\
&&+2\,\Delta t\left(\nor{P_{R}^{\prime}(\uv_{r}^{n+1}\cdot\nabla\boldsymbol{\eta}^{n+1})}{\tau}^{2}+\nor{P_{R}^{\prime}(\uv_{r}^{n+1}\cdot\nabla\uv^{n+1})}{\tau}^{2}\right),
\end{eqnarray}
where $C_{1}$ is a constant depending on $\nu^{-1}$.

\medskip

Summing \eqref{eq:ErrEq3} from $n=0$ to $k\leq N-1$, we have:
\begin{eqnarray}\label{eq:ErrEqSum}
&&\max_{0\leq k\leq N}\nor{\boldsymbol{\phi}_{r}^{k}}{{\bf L}^{2}}^{2}
+\nu\Delta t\sum_{n=0}^{N-1}\nor{\nabla\boldsymbol{\phi}_{r}^{n+1}}{{\bf L}^{2}}^{2}
+\Delta t\sum_{n=0}^{N-1}\nor{P_{R}^{\prime}(\uv_{r}^{n+1}\cdot\nabla\boldsymbol{\phi}_{r}^{n+1})}{\tau}^{2}\nonumber \\
&\leq&\nor{\boldsymbol{\phi}_{r}^{0}}{{\bf L}^{2}}^{2}+C_{1}\,\Delta t\sum_{n=0}^{N-1}\nor{\boldsymbol{r}^{n}}{{\bf H}^{-1}}^{2}
+C_{1}\,\Delta t\sum_{n=0}^{N-1}\left\| \frac{1}{\Delta t}(\boldsymbol{\eta}^{n+1}-\boldsymbol{\eta}^{n})\right\|_{{\bf L}^2}^{2}\nonumber \\
&&+C_{1}\,\Delta t\sum_{n=0}^{N-1}\left(\nor{\nabla\uv_{r}^{n+1}}{{\bf L}^{2}}^{2}
+\nor{\nabla\uv^{n+1}}{{\bf L}^{2}}^{2}\right)\nor{\nabla\boldsymbol{\eta}^{n+1}}{{\bf L}^{2}}^{2}\nonumber \\
&&+C_{1}\,\Delta t\sum_{n=0}^{N-1}\nor{\nabla\uv^{n+1}}{{\bf H}^{1}}^{2}\nor{\boldsymbol{\phi}_{r}^{n+1}}{{\bf L}^{2}}^{2} +C_{1}\,\Delta t \sum_{n=0}^{N-1}\nor{p^{n+1}-q_{h}}{{\bf L}^{2}}^{2}
\nonumber \\
&&+2\,\Delta t\sum_{n=0}^{N-1}\left(\nor{P_{R}^{\prime}(\uv_{r}^{n+1}\cdot\nabla\boldsymbol{\eta}^{n+1})}{\tau}^{2}+\nor{P_{R}^{\prime}(\uv_{r}^{n+1}\cdot\nabla\uv^{n+1})}{\tau}^{2}\right).
\end{eqnarray}
Next, we estimate each term on the right-hand side of \eqref{eq:ErrEqSum}.

\medskip

The first term on the right-hand side of \eqref{eq:ErrEqSum} can be estimated as follows:
\BEQ\label{eq:ErrIneq2RHS1}
\nor{\boldsymbol{\phi}_{r}^{0}}{}^{2}\leq \nor{\uv^{0}-\uv_{r}^{0}}{}^{2}+\nor{\uv^{0}-\vv_{r}^{0}}{}^{2}\leq C\,h^{2\ell},
\EEQ
where the last inequality follows from the fact that $\vv_{r}^{0}$ is the elliptic projection of $\uv^{0}$ on $\Xv^{r}\subset \Vv^{h}$, so that it satisfies optimal approximation properties similar to standard FE interpolations (cf. \cite{Ciarlet02}), and we have supposed $\nor{\uv^{0}-\uv_{r}^{0}}{}=\mathcal{O}(h^{\ell})$.

\medskip

By using Poincar\'e-Friedrichs inequality, the second term on the right-hand side of \eqref{eq:ErrEqSum} can be estimated as follows (see, e.g., \cite{IliescuWang13}):
\BEQ\label{eq:rn}
\Delta t\sum_{n=0}^{N-1}\nor{\boldsymbol{r}^{n}}{{\bf H}^{-1}}^{2}\leq C\,\Delta t\sum_{n=0}^{N-1}\nor{\boldsymbol{r}^{n}}{{\bf L}^{2}}^{2}\leq C\,\Delta t^{2}\nor{\partial_{t}^{2}\uv}{L^{2}({\bf L}^{2})}^{2}.
\EEQ

Using Corollary \ref{co:PODco}, the third term on the right-hand side of \eqref{eq:ErrEqSum} can be estimated as follows:
\BEQ\label{eq:eta1}
\Delta t\sum_{n=0}^{N-1}\left\| \frac{1}{\Delta t}(\boldsymbol{\eta}^{n+1}-\boldsymbol{\eta}^{n})\right\|_{{\bf L}^2}^{2}
\leq \nor{\partial_{t}\boldsymbol{\eta}}{L^{2}({\bf L}^{2})}^{2}\leq C\left(h^{2\ell} + \Delta t ^{2} + \sum_{i=r+1}^{M}\lambda_{i}\right).
\EEQ
To estimate the fourth term on the right-hand side of \eqref{eq:ErrEqSum}, we use Lemma \ref{lm:stabEst} and the fact that $\vv_{r}^{n+1}$ is the elliptic projection of $\uv^{n+1}$ on $\Xv^{r}\subset\Vv^{h}$, so that it satisfies optimal approximation properties as standard FE interpolations (cf. \cite{Ciarlet02}):
\BEQ\label{eq:eta2}
\Delta t\sum_{n=0}^{N-1}\nor{\nabla\uv_{r}^{n+1}}{{\bf L}^{2}}^{2}\nor{\nabla\boldsymbol{\eta}^{n+1}}{{\bf L}^{2}}^{2}
\leq C_{2}\,h^{2\ell},
\EEQ
where $C_2$ is a constant depending on $\nu^{-2}$.

\medskip

By using the regularity assumption \ref{hp:FERegularity} on the continuous solution and \eqref{eq:PODlmH01}, the fifth term on the right-hand side of \eqref{eq:ErrEqSum} can be estimated as follows:
\begin{eqnarray}\label{eq:eta3}
\Delta t\sum_{n=0}^{N-1}\nor{\nabla\uv^{n+1}}{{\bf L}^{2}}^{2}\nor{\nabla\boldsymbol{\eta}^{n+1}}{{\bf L}^{2}}^{2}
\leq C\left(h^{2\ell}+\Delta t ^{2} + \sum_{i=r+1}^{M}\lambda_{i}\right).
\end{eqnarray}
Since \eqref{eq:presErrEq} holds for any $q_{h}\in Q_{h}$, we can use the pressure approximation property \eqref{eq:FEpresEE} in Hypothesis \ref{hp:FEerrEst} to bound the seventh term on the right-hand side of \eqref{eq:ErrEqSum}:
\BEQ\label{eq:pn+1}
\Delta t \sum_{n=0}^{N-1}\nor{p^{n+1}-q_{h}}{{\bf L}^{2}}^{2}\leq C\,h^{2\ell}.
\EEQ
Using Lemma \ref{lm:stabPR}, Minkowski's and H\"older's inequalities, the eighth term on the right-hand side of \eqref{eq:ErrEqSum} can be estimated as follows:
\begin{eqnarray*}
&&\Delta t\sum_{n=0}^{N-1}\nor{P_{R}^{\prime}(\uv_{r}^{n+1}\cdot\nabla\boldsymbol{\eta}^{n+1})}{\tau}^{2}
\leq C\,\Delta t\,h^{2}\,\sum_{n=0}^{N-1}\nor{(\uv_{r}^{n+1}\pm\vv_{r}^{n+1})\cdot\nabla\boldsymbol{\eta}^{n+1}}{{\bf L}^{2}}^{2} \\
&\leq& C\,\Delta t\,h^{2}\,\sum_{n=0}^{N-1}\left(\nor{\boldsymbol{\phi}_{r}^{n+1}\cdot\nabla\boldsymbol{\eta}^{n+1}}{{\bf L}^{2}}^{2}+\nor{\vv_{r}^{n+1}\cdot\nabla\boldsymbol{\eta}^{n+1}}{{\bf L}^{2}}^{2}\right) \\
&\leq& C\,\Delta t\,h^{2}\,\sum_{n=0}^{N-1}\left(\nor{\boldsymbol{\phi}_{r}^{n+1}}{{\bf L}^{6}}^{2}\nor{\nabla\boldsymbol{\eta}^{n+1}}{{\bf L}^{3}}^{2}+\nor{\vv_{r}^{n+1}\pm\uv^{n+1}}{{\bf L}^{6}}^{2}\nor{\nabla\boldsymbol{\eta}^{n+1}}{{\bf L}^{3}}^{2}\right) \\
&\leq& C\,\Delta t\,h^{2}\,\sum_{n=0}^{N-1}\left(\nor{\nabla\boldsymbol{\phi}_{r}^{n+1}}{{\bf L}^{2}}^{2}\nor{\nabla\boldsymbol{\eta}^{n+1}}{{\bf L}^{3}}^{2}+
\left(\nor{\nabla\uv^{n+1}}{{\bf L}^{2}}^{2}+\nor{\nabla\boldsymbol{\eta}^{n+1}}{{\bf L}^{2}}^{2}\right)\nor{\nabla\boldsymbol{\eta}^{n+1}}{{\bf L}^{3}}^{2}\right),
\end{eqnarray*}
where the last inequality comes from the Sobolev embedding ${\bf H}^{1}\hookrightarrow {\bf L}^{6}$. Now, using optimal approximation properties for standard FE interpolations (cf. \cite{Ciarlet02}) and local inverse estimates (cf. \cite{Bernardi04}), we have that $\nor{\nabla\boldsymbol{\eta}^{n+1}}{{\bf L}^{3}}^{2}\leq C\left(h^{2\ell -1}+h^{-1}\nor{\nabla\boldsymbol{\eta}^{n+1}}{{\bf L}^{2}}^{2}\right)$. Using this fact in the above inequality, we get:
\begin{eqnarray*}
&&\Delta t\sum_{n=0}^{N-1}\nor{P_{R}^{\prime}(\uv_{r}^{n+1}\cdot\nabla\boldsymbol{\eta}^{n+1})}{\tau}^{2}\\
&\leq&
C\,\Delta t\,\sum_{n=0}^{N-1}\left(\nor{\nabla\boldsymbol{\phi}_{r}^{n+1}}{{\bf L}^{2}}^{2}
+\nor{\nabla\uv^{n+1}}{{\bf L}^{2}}^{2}+\nor{\nabla\boldsymbol{\eta}^{n+1}}{{\bf L}^{2}}^{2}\right)
\left(h^{2\ell +1}+h\nor{\nabla\boldsymbol{\eta}^{n+1}}{{\bf L}^{2}}^{2}\right).
\end{eqnarray*}
Taking $h\sim \Delta t$, using \eqref{eq:PODlmH01} and the regularity assumption \ref{hp:FERegularity}, we obtain:
\begin{eqnarray}\label{eq:stabErrEst1}
&&\Delta t\sum_{n=0}^{N-1}\nor{P_{R}^{\prime}(\uv_{r}^{n+1}\cdot\nabla\boldsymbol{\eta}^{n+1})}{\tau}^{2}\nonumber \\
&\leq&
C\,\Delta t\,\sum_{n=0}^{N-1}\nor{\nabla\boldsymbol{\phi}_{r}^{n+1}}{{\bf L}^{2}}^{2}
\left(h^{2\ell}+\Delta t^{2} + \sum_{i=r+1}^{M}\lambda_{i}\right)\nonumber \\
&&+C\, \left(h^{2\ell}+\Delta t^{2} + \sum_{i=r+1}^{M}\lambda_{i}\right).
\end{eqnarray}
Using again Lemma \ref{lm:stabPR}, Minkowski's and H\"older's inequalities, we have the following error bound for the last term on the right-hand side of \eqref{eq:ErrEqSum}:
\begin{eqnarray*}
&&\Delta t\sum_{n=0}^{N-1}\nor{P_{R}^{\prime}\left((\uv_{r}^{n+1}\pm\vv_{r}^{n+1})\cdot\nabla\uv^{n+1}\right)}{\tau}^{2}\\
&\leq& C\,\Delta t\,\sum_{n=0}^{N-1}h^{2}\nor{\boldsymbol{\phi}_{r}^{n+1}\cdot\nabla\uv^{n+1}}{{\bf L}^{2}}^{2}
+\nor{P_{R}^{\prime}\left((\vv_{r}^{n+1}\pm\uv^{n+1})\cdot\nabla\uv^{n+1}\right)}{\tau}^{2}\\
&\leq& C\,\Delta t\,h^{2}\,\sum_{n=0}^{N-1}\left(\nor{\boldsymbol{\phi}_{r}^{n+1}}{{\bf L}^{6}}^{2}\nor{\nabla\uv^{n+1}}{{\bf L}^{3}}^{2}
+\nor{\boldsymbol{\eta}^{n+1}}{{\bf L}^{6}}^{2}\nor{\nabla\uv^{n+1}}{{\bf L}^{3}}^{2}\right)\\
&&+\Delta t\,\sum_{n=0}^{N-1}\nor{P_{R}^{\prime}(\uv^{n+1}\cdot\nabla\uv^{n+1})}{\tau}^{2}\\
&\leq& C\,\Delta t\,h^{2}\,\sum_{n=0}^{N-1}\left(\nor{\nabla\boldsymbol{\phi}_{r}^{n+1}}{{\bf L}^{2}}^{2}
+\nor{\nabla\boldsymbol{\eta}^{n+1}}{{\bf L}^{2}}^{2}\right)+\Delta t\,\sum_{n=0}^{N-1}\nor{P_{R}^{\prime}(\uv^{n+1}\cdot\nabla\uv^{n+1})}{\tau}^{2},
\end{eqnarray*}
where the last inequality comes from the Sobolev embedding ${\bf H}^{1}\hookrightarrow {\bf L}^{6}$, and the regularity assumption \ref{hp:FERegularity} together with the Sobolev embedding ${\bf H}^{2}\hookrightarrow {\bf W}^{1,3}$. Now, using \eqref{eq:PODlmH01}, we get:
\begin{eqnarray*}
&&\Delta t\sum_{n=0}^{N-1}\nor{P_{R}^{\prime}\left(\uv_{r}^{n+1}\cdot\nabla\uv^{n+1}\right)}{\tau}^{2}\\
&\leq& C\,\Delta t\,h^{2}\,\sum_{n=0}^{N-1}\nor{\nabla\boldsymbol{\phi}_{r}^{n+1}}{{\bf L}^{2}}^{2}
+C\, \left(h^{2\ell}+\Delta t^{2} + \sum_{i=r+1}^{M}\lambda_{i}\right)
\\
&&+\Delta t\,\sum_{n=0}^{N-1}\nor{P_{R}^{\prime}(\uv^{n+1}\cdot\nabla\uv^{n+1})}{\tau}^{2}.
\end{eqnarray*}
Here, the last term can be bounded as follows:
\begin{eqnarray*}
\Delta t\,\sum_{n=0}^{N-1}\nor{P_{R}^{\prime}(\uv^{n+1}\cdot\nabla\uv^{n+1})}{\tau}^{2}
&=&\Delta t\,\sum_{n=0}^{N-1}\nor{(\uv^{n+1}\cdot\nabla\uv^{n+1})-P_{R}(\uv^{n+1}\cdot\nabla\uv^{n+1})}{\tau}^{2}\\
&\leq& C\left(h^{2l}+\Delta t ^{2} + h^{2}\sum_{i=R+1}^{M}\widehat{\lambda}_{i}\right),
\end{eqnarray*}
where we have used assumption \eqref{eq:stabCoef} in Hypothesis \ref{hp:estimCoef} on the stabilization parameters, assumption \eqref{eq:FEvelEE} in Hypothesis \ref{hp:FEerrEst} on the FE velocity approximation, \eqref{eq:PR}, \eqref{eq:PODlmL2}, and $\widehat{\lambda}_{i}$, $i=R+1,\ldots,M$ are the eigenvalues associated to the snapshots correlation matrix $\widehat{K}$ previously defined in \eqref{eq:Khat}.
Thus, the last term on the right-hand side of \eqref{eq:ErrEqSum} is finally bounded as:
\begin{eqnarray}\label{eq:stabErrEst2}
&&\Delta t\sum_{n=0}^{N-1}\nor{P_{R}^{\prime}\left(\uv_{r}^{n+1}\cdot\nabla\uv^{n+1}\right)}{\tau}^{2}\nonumber \\
&\leq &C\,\Delta t\,h^{2}\,\sum_{n=0}^{N-1}\nor{\nabla\boldsymbol{\phi}_{r}^{n+1}}{{\bf L}^{2}}^{2}\nonumber \\
&&+C\, \left(h^{2\ell}+\Delta t^{2} + \sum_{i=r+1}^{M}\lambda_{i}+h^{2}\sum_{i=R+1}^{M}\widehat{\lambda}_{i}\right).
\end{eqnarray}
Collecting \eqref{eq:rn}-\eqref{eq:stabErrEst2}, by dropping the third term on the left-hand side of \eqref{eq:ErrEqSum}, this latter becomes:
\begin{eqnarray}\label{eq:ErrEqRiep}
&&\max_{0\leq k\leq N}\nor{\boldsymbol{\phi}_{r}^{k}}{{\bf L}^{2}}^{2} +\Delta t\left[\nu-C(h^{2} + E)\right]\sum_{n=0}^{N-1}\nor{\nabla\boldsymbol{\phi}_{r}^{n+1}}{{\bf L}^{2}}^{2}
\nonumber \\ 
&\leq&C_{1}\,\Delta t\sum_{n=0}^{N-1}\nor{\nabla\uv^{n+1}}{{\bf H}^{1}}^{2}\nor{\boldsymbol{\phi}_{r}^{n+1}}{{\bf L}^{2}}^{2}
+C_{3}\left( E+h^{2}\sum_{i=R+1}^{M}\widehat{\lambda}_{i}\right),
\end{eqnarray}
where $C_{3}$ is a constant depending on $\nu^{-3}$, and we have called:
\BEQ\label{eq:E}
E=\left(h^{2\ell}+\Delta t^{2} + \sum_{i=r+1}^{M}\lambda_{i}\right).
\EEQ
For a sufficiently small $h\sim \Delta t\left(\sim\mathcal{O}(\sqrt{\nu})\right)$ and a sufficiently large $r:\lambda_{r+1}\sim\mathcal{O}(\nu)$ , $\nu-C(h^{2} + E)\geq \nu/2$. If $\Delta t\leq \Delta t^{*}=\disp\frac{1}{2C_{1}\disp\max_{n=0,\ldots,N-1}\left(\nor{\nabla\uv^{n+1}}{{\bf H}^{1}}^{2}\right)}$, the discrete Gr\"onwall's lemma (see Lemma 27 in \cite{Layton08} for instance) implies the following inequality:
\BEQ\label{eq:ErrEqRiepF}
\max_{0\leq k\leq N}\nor{\boldsymbol{\phi}_{r}^{k}}{{\bf L}^{2}}^{2}
+\frac{\nu}{2}\Delta t\sum_{n=0}^{N-1}\nor{\nabla\boldsymbol{\phi}_{r}^{n+1}}{{\bf L}^{2}}^{2}
\leq C^{*}\left(E+h^{2}\sum_{i=R+1}^{M}\widehat{\lambda}_{i}\right),
\EEQ
where $C^{*}=C_{3}e^{C_{1}\Delta t\sum_{n=0}^{N-1}\nor{\nabla\uv^{n+1}}{{\bf H}^{1}}^{2}}$. Finally, using in \eqref{eq:ErrEqRiepF} the obvious inequality:
$$
\max_{0\leq k\leq N}\nor{\boldsymbol{\phi}_{r}^{k}}{{\bf L}^{2}}^{2}\geq
\frac{1}{N+1}\sum_{n=0}^{N}\nor{\boldsymbol{\phi}_{r}^{n}}{{\bf L}^{2}}^{2},
$$
triangle inequality and estimates \eqref{eq:PODlmL2}-\eqref{eq:PODlmH01}, we get:
\BEQ\label{eq:PODEEF}
\frac{1}{N+1}\sum_{n=0}^{N}\nor{\uv^{n}-\uv_{r}^{n}}{{\bf L}^{2}}^{2}+\nu\Delta t\sum_{n=0}^{N-1}\nor{\nabla(\uv^{n+1}-\uv_{r}^{n+1})}{{\bf L}^{2}}^{2}
\leq C\left(E + \sum_{i=r+1}^{M}\lambda_{i}+h^{2}\sum_{i=R+1}^{M}\widehat{\lambda}_{i}\right).
\EEQ
This concludes the proof.
\qed
\begin{theorem}\label{th:SI}
Under the hypotheses of Theorem \ref{th:PODEE}, the solution of the semi-implicit SD-POD-ROM \eqref{eq:SD-POD-ROMdiscSI} satisfies error estimate \eqref{eq:PODEE}.
\end{theorem}
The proof of this theorem can be achieved by the same techniques used to prove Theorem \ref{th:PODEE}, thus we skip it for brevity.
\begin{remark}\label{rm:PS-POD-ROM2}
If one consider $\mathcal{H}={\bf L}^{2}$ in the generation of POD modes (as in \cite{IliescuWang14} for instance), then the following error estimate can be derived for the SD-POD-ROM \eqref{eq:SD-POD-ROMdisc} or \eqref{eq:SD-POD-ROMdiscSI}:
\begin{eqnarray}\label{eq:PODEEl2}
&&\frac{1}{N+1}\sum_{n=0}^{N}\nor{\uv^{n}-\uv_{r}^{n}}{{\bf L}^{2}}^{2}+\nu\Delta t\sum_{n=0}^{N-1}\nor{\nabla(\uv^{n+1}-\uv_{r}^{n+1})}{{\bf L}^{2}}^{2}\nonumber \\
&\leq& C\left((1+\nor{S_{r}}{2})(h^{2\ell}+\Delta t^{2}) + \sum_{i=r+1}^{M}(1+\nor{\nabla\boldsymbol{\varphi}_{i}}{{\bf L}^{2}}^{2})\lambda_{i}+h^{2}\sum_{i=R+1}^{M}\widehat{\lambda}_{i}\right),
\end{eqnarray}
with $\nor{S_{r}}{2}$ denoting the $2$-norm of the stiffness matrix with entries $[S_{r}]_{ij}=(\nabla\boldsymbol{\varphi}_{j},\nabla\boldsymbol{\varphi}_{i})$, $i,j=1,\ldots,r$. The appearance of $\nor{S_{r}}{2}$ comes from the use of the POD inverse estimate (see \cite{KunischVolkwein01}, Lemma 2):
\BEQ\label{eq:PODinvest}
\nor{\nabla\vv}{{\bf L}^{2}}\leq \sqrt{\nor{S_{r}}{2}}\nor{\vv}{{\bf L}^{2}},\quad \forall \vv\in\Xv^{r}.
\EEQ
In this case, there is no need to include the finite difference quotients in the set of snapshots to prove optimal error estimates in time as for the case $\mathcal{H}={\bf H}_{0}^{1}$. We notice, however, that for practical computations, one would use rather fine time discretizations for snapshots, for which the inclusion of the difference quotients in the case $\mathcal{H}={\bf H}_{0}^{1}$ should be almost unnoticeable (see numerical evidences in \cite{KunischVolkwein01} for instance).
\end{remark}

\begin{remark}\label{rm:PS-POD-ROM2}
If one just consider the standard Galerkin POD-ROM ($\tau_{K}=0$ for any $K\in {\cal T}_{h}$), thus error estimate \eqref{eq:PODEE} can be recovered, without the appearance of the last term on the right-hand side of \eqref{eq:PODEE}. In this case, any control on the high-frequency modes of the convective derivative is guaranteed. 
When $R=0$, one has that the last term on the right-hand side of \eqref{eq:PODEE} is limited to $h^2$. In this case, there is no interest in increasing more than $\ell=2$ the order of the FE velocity interpolation to construct the POD basis, as this would requires a larger computational effort without increasing the accuracy of the POD-ROM numerical solution. This low convergence order appears linked to the diffusive nature of the penalty-stabilized POD-ROM \eqref{eq:PS-POD-ROM}, which extends the numerical diffusion to all the resolved modes.
\end{remark}

\begin{remark}
Note that to prove estimate \eqref{eq:PODEE}, we have to assume that  \eqref{eq:FEvelEE} in Hypothesis \ref{hp:FEerrEst} holds for the FE velocity approximation. This optimal convergence order assumption is generally valid in laminar flow settings or for sufficiently regular flows, but is usually not valid in realistic turbulent flow settings, since the convergence order decreases with the regularity of the flow. Finding robust numerical schemes for realistic turbulent flows is still an open issue, to the best of our knowledge ({cf.} \cite{Berselli06, Layton08}). However, the main goal of this report is not to develop robust numerical schemes for turbulent flows. As pointed out also in \cite{IliescuWang14}, we assume that an acceptable scheme exists and we investigate whether
the stabilized POD-ROM that we are considering, with features particularly sui\-table in the convection-dominated regime, can achieve a similar numerical accuracy, but with a dramatically reduced dimension. This is a common approach in the derivation of error estimates for POD-ROM ({\em cf.} \cite{IliescuWang14, Luo08}). The achievement of an optimal numerical accuracy is subject to the condition $r:\lambda_{r+1}\sim\mathcal{O}(\nu)$, which gives an idea on how many POD modes are needed to reach this accuracy for a certain fluid viscosity $\nu$. 
\end{remark}

\section{Practical implementation}\label{sec:PI}

In this section, we suggest an efficient practical implementation of the stabilization term, where the stabilization parameter is approximated by the Discrete Empirical Interpolation Method (DEIM, {\em cf.} \cite{Chaturantabut10}). This leads to a computationally efficient and mathematically founded offline/online algorithm (completely separated), implemented over the standard POD-Galerkin ROM. To the best of the authors' knowledge, the use of DEIM for the accurate and efficient computation of the stabilization parameter is new in the literature so far. Indeed, in \cite{IliescuJohn15} a stabilization parameter simply arising from the FE resolution is compared towards a stabilization parameter just based on the POD spatial resolution for a SUPG-ROM applied to advection-diffusion-reaction equations, while in \cite{BergmannIollo09} an optimization problem for the determination of the stabilization parameter is solved.

\medskip

The proposed strategy consists in approximating the local stabilization parameters as a piecewise constant FE function $\tau$ reading as:
\BEQ\label{eq:tau}
\tau=\sum_{m=1}^{\widetilde{r}}\alpha_{m}(t)\rho_{m}(\xv),
\EEQ
with $\widetilde{r}\sim \mathcal{O}(r)$, and $\alpha_{m}(t), \rho_{m}(\xv)$ that will be determined in the next section.
 
\subsection{DEIM algorithm for the computation of the stabilization parameter}\label{sec:POD-ROM}

To describe the DEIM algorithm for the computation of the stabilization parameter $\tau$ in \eqref{eq:tau}, we are going to follow the notation used in \cite{Quarteroni16}, Section 10.3.
\\
\begin{itemize}
\item{OFFLINE PHASE.}
\end{itemize}

\begin{itemize}
\item[(i)]
The offline phase consists first in constructing the spatial basis $\mathbb{Q}=\left[\rho_{1}|\ldots|\rho_{\widetilde{r}}\right]$, obtained by operating a Singular Value Decomposition (SVD) on a set of snapshots $\left[\tau_{h}(\cdot,t_{1}),\ldots,\tau_{h}(\cdot,t_{N}))\right]$, where we may use the following expression for the offline piecewise constant FE stabilization coefficient:
\BEQ\label{eq:tauh}
\tau_{h}(\cdot,t_{i})=\left[c_{1}\frac{\nu}{h_{K}^{2}}+c_{2}\frac{U_{K}^{i}}{h_{K}}\right]^{-1},\quad i=1,\ldots,N.
\EEQ
In \eqref{eq:tauh}, $c_{1}$ and $c_{2}$ are user-chosen positive constants, and $U_{K}^{i}, i=1,\ldots,N$, is some local speed on the mesh cell $K$ at the offline time step $t_{i}$, e.g.: 
$$
U_{K}^{i}=\nor{\uhv(t_{i})}{{\bf L}^{2}(K)}/|K|^{1/2}.
$$
Form \eqref{eq:tauh} is designed by asymptotic scaling arguments applied in the framework of stabilized methods ({\em cf.} \cite{Codina07}), aimed at taking into account the local balance between convection and diffusion. This ensures a self-adapting high accuracy up to high Reynolds number flows.
\item[(ii)]
The second step of the offline phase consists in selecting iteratively $\widetilde{r}$ indices $\mathcal{I}\subset\{1,\ldots,N_{K}\}$, where $N_{K}$ is the number of elements $K\in{\cal T}_{h}$, from the basis $\mathbb{Q}$ using the following greedy procedure, which minimizes at each step the interpolation error over the snapshots set measured in the maximum norm:
\begin{itemize}
\item {\bf Initialization:} $i_{1}=\text{arg}\disp\max_{i=1,\ldots,N_{K}} |(\rho_{1})_{i}|$; $\mathbb{Q}=\rho_{1}$; $\mathcal{I}=\{i_{1}\}$.
\item {\bf Iterations:} 
$$
\left\{
\begin{array}{ll}
\mbox{\bf for} & m=2:\widetilde{r}\\
& \mbox{res}=\rho_{m}-\mathbb{Q}\,\mathbb{Q}_{\mathcal{I}}^{-1}(\rho_{m})_{\mathcal{I}};\\
& i_{m}=\text{arg}\disp\max_{i=1,\ldots,N_{K}} |\mbox{res}_{i}|;\\
& \mathbb{Q}\leftarrow [\mathbb{Q}|\rho_{m}]; \mathcal{I}\leftarrow \mathcal{I}\cup\{i_{m}\};\\
\mbox{\bf end} & 
\end{array}
\right.
$$
where $\mathbb{Q}_{\mathcal{I}}$ is the matrix formed by the $\mathcal{I}$ rows of $\mathbb{Q}$, and $(\rho_{m})_{\mathcal{I}}$ is the piecewise constant FE vector formed by the $\mathcal{I}$ components of $\rho_{m}$.
\end{itemize}
\end{itemize}

\begin{itemize}
\item{ONLINE PHASE.}
\end{itemize}

\begin{itemize}
\item[(i)]
In order to compute online the time coefficients vector $\boldsymbol{\alpha}(t)=[\alpha_{1}(t),\ldots,\alpha_{\widetilde{r}}(t)]^{T}$, interpolation constraints are imposed at the $\widetilde{r}$ points corresponding to the selected indices. So, first we form ${\boldsymbol{\tau}}_{\mathcal{I}}^{n}$ by evaluating: 
\BEQ\label{eq:taur}
\tau_{r}(\cdot,t_{n})=\left[c_{1}\frac{\nu}{h_{K}^{2}}+c_{2}\frac{\nor{\uv_{r}^{n}}{{\bf L}^{2}(K)}/|K|^{1/2}}{h_{K}}\right]^{-1},
\EEQ
on the interpolation points $\{\xv^{i_{1}},\ldots,\xv^{i_{\widetilde{r}}}\}$. Note that online time steps may differ from the offline ones (remember that the proposed method is a stabilized POD-ROM with time as parameter).
\item[(ii)]
The second step of the online phase requires the solution of the following linear system:
\BEQ\label{eq:DEIM}
\mathbb{Q}\boldsymbol{\alpha}^{n}={\boldsymbol{\tau}}_{\mathcal{I}}^{n},
\EEQ
which gives the searched time coefficients vector $\boldsymbol{\alpha}^{n}=[\alpha_{1}^{n},\ldots,\alpha_{\widetilde{r}}^{n}]^{T}$, and has comple\-xity $\mathcal{O}(\widetilde{r}^{\,3})$.
\end{itemize}

\section{Summary and conclusions}\label{sec:Concl}

In this work, we have proposed a new stabilized POD-ROM for the numerical simulation of convection-dominated incompressible fluid flows. This model, denoted SD-POD-ROM, is derived from the VMS formalism, and uses a streamline derivative projection-based ope\-rator to properly take into account the high frequencies convective derivative component of POD modes not included in the ROM.

\medskip

We have performed a stability and convergence analysis of the arising fully discrete SD-POD-ROM applied to the unsteady incompressible NSE. The main contribution of the present paper is the proof of a sharp error estimate that considers all contributions: the spatial discretization error (due to the FE discretization), the temporal discretization error (due to the backward Euler method), and the POD truncation error. In particular, the numerical analysis makes apparent an extra-control on the high frequencies of the convective derivative, which is an extremely important feature in view of computing turbulent flows. The question of an efficient practical implementation of the strongly non-linear convective stabilization term within the SD-POD-ROM is also addressed, using DEIM to approximate the non-linear stabilization parameter.

\medskip

We plan to extend this theoretical work on the numerical analysis of the proposed SD-POD-ROM by performing a numerical investigation that both supports the analytical results and illustrate the potential of the method for the challenging simulation of turbulent flows. This computational study is today in progress, and shall appear in a forthcoming paper.

\medskip

{\sl Acknowledgments:} The research of Tom\'as Chac\'on Rebollo and Samuele Rubino has been partially funded by the Spanish Government - EU FEDER Project MTM2015-64577-C2-1-R. Samuele Rubino would also gratefully acknowledge the financial support received from IdEx (Initiative d'Excellence de l'Universit\'e de Bordeaux) International Post-Doc Program during his postdoctoral research involved in this article.

\bibliographystyle{abbrv}
\bibliography{Biblio_SD-POD-ROM}


\end{document}